\documentclass[reqno]{article}

\def\p{\varphi}

\def\D{\partial}
\def\l{\ell}

\def\ep{\epsilon}
\def\a{\alpha}

\def\d{\delta}
\def\be{\begin{equation}}
\def\ee{\end{equation}}
\def\bea{\begin{eqnarray}}
\def\eea{\end{eqnarray}}
\def\ba{\begin{array}}
\def\ea{\end{array}}

\def\bx{{\bf x}}

\def\bc{{\bf c}}

\def\bxi{{\bf \xi}}

\def\ra{\rightarrow}
\def\R{{\mathbb R}}

\def\C{{\mathbb C}}

\usepackage{amsmath}
\usepackage{amsfonts}
\usepackage{amsthm}

\newtheorem{theorem}{Theorem}[section]
\newtheorem{lemma}[theorem]{Lemma}

\theoremstyle{definition}

\theoremstyle{remark}

\numberwithin{equation}{section}

\newcommand{\abs}[1]{\lvert#1\rvert}

\newfont{\Bb}{msbm8 scaled\magstep{1}}
\newcommand{\rc}{\mbox{\Bb R}}

\begin{document}

\title{Modified Rayleigh Conjecture Method and Its Applications 
}

\author{
Alexander G. RAMM\footnote{Department of Mathematics, Kansas State University,
Manhattan, Kansas 66506-2602, USA, e-mail: ramm@math.ksu.edu}\\
and
\\
Semion GUTMAN\footnote{Department of Mathematics, University of Oklahoma, Norman,
 OK 73019, USA, e-mail: sgutman@ou.edu}
}
\date{}
\maketitle

\begin{abstract} 
 The Rayleigh conjecture about convergence up to the boundary of the 
series representing 
the scattered 
field in the exterior of an obstacle $D$  is widely used by engineers in 
applications.
However this conjecture is false for some obstacles. 
AGR introduced the Modified Rayleigh
Conjecture (MRC), which is an exact mathematical result.
In this paper we review the theoretical basis for the MRC method for 2D and 3D
obstacle scattering problems, for static problems, and for scattering by
periodic structures. We also present successful numerical
algorithms based on the MRC for various scattering problems.
The MRC method is easy
to implement for both simple and complex geometries. It is shown to be a
viable alternative for other obstacle scattering methods. Various direct
and inverse scattering problems require finding global minima of
functions of several variables. The Stability Index Method (SIM) combines
stochastic and deterministic method to accomplish such a minimization.

\end{abstract}

\vskip.3cm

{\bf Math. Subj. classification:} 35J05, 65M99, 78A40





\section{Introduction}

In this paper we review our recent results on the Modified Rayleigh
Conjecture (MRC) method. The method is applied to multidimensional
obstacle scattering problems, as well as to scattering by periodic
structures. Also we discuss an application of the MRC to static
problems, and preliminary results on inverse obstacle scattering by MRC.
Numerical results illustrate the performance of various MRC algorithms.
The paper concludes with a presentation of the Stability Index
Method (SIM) for global minimization.

The basic theoretical foundation of the
method was developed in \cite{r430}. The MRC has the appeal
of an easy implementation for obstacles of complicated
geometry, e.g. having edges and corners. In our numerical
experiments the method has shown itself to be a competitive
alternative to the BIEM (boundary integral equations
method), see \cite{r437}. Also, unlike the BIEM, one can
apply the algorithm to different obstacles with very little
additional effort. A similar method is discussed in \cite{doicu}
 
We formulate the obstacle scattering problem in a 3D setting
with the
 Dirichlet boundary condition, but the method can
also be used for the Neumann boundary 
condition, corresponding to acoustically hard obstacles, and
the Robin  boundary condition.

Consider a bounded domain $D \subset \R^3$, with a
 Lipschitz boundary $S$.
Denote the
exterior domain by $D^\prime = \R^3 \backslash D$. Let
 $\alpha, \alpha^\prime\in S^2$ be unit vectors,
where $ S^2$ is the unit sphere in $\R^3$.

The acoustic wave scattering problem by an acoustically soft obstacle $D$ 
consists in
finding the (unique) solution to the problem \eqref{s1_e1}-\eqref{s1_e2}:
\be\label{s1_e1}
\left(\nabla^2 + k^2 \right) u=0 \hbox{\ in\ } D^\prime, \quad
  u = 0 \hbox{\ on\ } S, 
\end{equation}

\be\label{s1_e2}
u=u_0 + A(\alpha^\prime, \alpha) \frac{e^{ikr}}{r} 
+ o
  \left(\frac{1}{r} \right), \quad r:=|x| \to \infty, \quad
  \alpha^\prime := \frac{x}{r}.
\end{equation}
Here $u_0:=e^{ik \alpha \cdot x}$ is the incident field,
$v:=u-u_0$ is the scattered field,
$A(\alpha^\prime, \alpha)$ is called the scattering amplitude, its
$k$-dependence is not shown, $k>0$ is the wavenumber. The scattered field
$v$ is an outgoing solution of the Helmholtz differential 
equation
\eqref{s1_e1}, that is, a solution which satisfies the
 radiation condition
\be\label{s1_somm}
\lim_{r\to\infty}\int_{|x|=r}\bigg|\frac{\partial v}{\partial
\abs{x}}-ikv\bigg|^2ds=0\,.
\end{equation}

Denote
\be\label{s1_e3}
A_\l (\alpha) := \int_{S^2} A(\alpha^\prime, \alpha)
  \overline{Y_\l (\alpha^\prime)} d\alpha^\prime, 
\end{equation}
where $Y_\l (\alpha)$ are the orthonormal spherical harmonics,
$Y_\l = Y_{\l m}, -\l \leq m \leq \l$.

Let a ball
$B_R := \{x : |x| \leq R\}$ contain the obstacle $D$.
Let $h_\l (r)$ be the
spherical Hankel functions, normalized so that $h_\l (r) \sim
\frac{e^{ikr}}{r}$ as $r \to +\infty$.
In the region $r>R$ the solution to (\ref{s1_e1})-(\ref{s1_e2}) is:
\be\label{ss2_e1}
u(x, \alpha) = e^{ik\alpha \cdot x} + 
\sum^\infty_{\l =0} A_\l (\alpha)\Psi_\l(x), \quad
\Psi_\l(x):= Y_\l (\alpha^\prime) h_\l (kr), \quad r > R,\quad  
\alpha^\prime =
  \frac{x}{r}, 
\end{equation}
where $r=|x|$, the sum includes the summation with respect to $m$, $-\l \leq m \leq \l$,
and $A_\l (\alpha)$ are defined in (\ref{s1_e3}), 
see \cite{r313}.

{\bf The Rayleigh conjecture (RC)} is: {\it the series
(\ref{ss2_e1}) converges up to the boundary $S$} (originally
RC dealt with periodic structures, gratings). This
conjecture is false for many obstacles, but is true for some
(\cite{baran,millar1,rammb1}). For example, if $n=2$ and $D$
is an ellipse, then the series analogous to (\ref{ss2_e1})
converges in the region $|x| >a$, where $2a$ is the distance
between the foci of the ellipse \cite{baran}. In the
engineering literature there are numerical algorithms based
on the Rayleigh conjecture.  These algorithms use projection
methods and are reported to be unstable. Moreover, no error
estimate has been obtained for such algorithms.  These
algorithms cannot converge for arbitrary obstacles, because
the Rayleigh conjecture is false for some obstacles.

Our aim is to give a formulation of a {\it Modified Rayleigh
Conjecture} (MRC) which holds for any Lipschitz obstacle and
can be used in numerical solution of direct and inverse
scattering problems. In other words, while the MRC still has
the word "conjecture" in its name, it is a proven
mathematical result for the scattered field in the exterior
domain $D'$. In contrast to algorithms based on the invalid
Rayleigh Conjecture,
 the MRC-based algorithms, like the ones described here,
converge and an error estimate for the approximate solution
they yield has been obtained in \cite{r430} (see also
\cite{r476}, Chapter 12). This 
error
estimate is sharp in the order $\epsilon$.

\section{Modified Rayleigh conjecture} What we call the
Modified Rayleigh Conjecture (MRC) is actually the following
Theorems \ref{s2_thm1}, see \cite{r430}, and \ref{s2_thm2},
see \cite{r493}. We denote by $H^m_{loc}(D')$ the set of
functions from the Sobolev space $H^m(\tilde {D})$ for any
compact strictly inner subdomain $\tilde{D}$ of $D'$, so
that the distance from $\tilde{D}$ to $S$ is positive, $dist 
(\tilde{D}, S)>0$.

\begin{theorem}\label{s2_thm1}

Let $v=u-u_0$ be the scattered field, where $u$ is 
the solution to (\ref{s1_e1})-(\ref{s1_e2}).
Then there exists a positive integer
$L=L(\ep)$ and the coefficients $c_\l=c_\l(\ep),\ 0\leq\l\leq L(\ep)$ such that

\begin{itemize}
   \item [(i).]
\be\label{s2_t1}
||u_0+v_\ep||_{L^2(S)} \leq \epsilon,
\end{equation}
where
\be\label{s2_t2}
v_\ep(x)=\sum_{\l=0}^{L(\epsilon)}c_\l(\epsilon)\Psi_\l(x).
\end{equation}

   \item [(ii).]
\be\label{s2_t3}
\|v_\ep-v\|_{L^2(S)}\leq\ep
\end{equation}
and
\be\label{s2_t4}
|||v_\ep-v|||=O(\ep)\,,\;\ep\ra 0\,,
\end{equation}
where
\[
|||\cdot|||=\|\cdot\|_{H^m_{loc}(D')}+\|\cdot\|_{L^2(D';(1+\abs{x})^{-
\gamma})}\,,
\]
$\gamma>1\,,\;m>0$ is an arbitrary integer.

   \item [(iii).]
\[
c_\l(\ep)\ra A_\l\,,\; \text{as}\  \ep\ra 0\,,\;  \forall \l \,,
\]
where $A_\l:=A_\l(\alpha)$ is defined in \eqref{s1_e3}.

 \end{itemize}
\end{theorem}

\begin{proof} First, we prove item (i). Then we establish Lemma
\ref{s2_lm1}, and continue with the proof of (ii) and (iii).

(i) Without loss of generality we can assume that the origin
is an interior point of the domain $D$.
To
establish \eqref{s2_t1} it is sufficient to show that
\be\label{s2_e63}
H:=\overline{span}\{\Psi_\l(s)\ :\ 0\leq\l<\infty,
\quad s\in S\}=L^2(S)\,.
\end{equation}
Suppose that there exists $p\in L^2(S),\; p\not=0$, such that $p\perp H$ 
in $L^2(S)$.
Define the single-layer potential by
\be\label{s2_e64}
W(y)=\int_S \frac{e^{ik|s-y|}}{ |s-y|} p(s)\ ds\,,\quad 
y\in\R^3,
\end{equation}
where $ds$ is the surface area element.
Let $U\subset D$ be a ball centered in the origin. 
Then the addition theorem for the fundamental solution implies that
$W(y)=0$ for any $y\in U$.

By the unique continuation
principle  $W\equiv 0$ in $D$. In particular $W=0$ on the boundary $S$.
Since $W$ is an outgoing solution of
$(\nabla^2 + k^2) W=0\; \text{in}\ D^\prime$ with $W = 0$ on $S$, one
concludes from the uniqueness of solutions to the Dirichlet
problem in $D'$ that $W\equiv 0$ in $\R^3$. Finally, the jump properties
of the
normal derivative of the single-layer potential imply 
that $p=0$ in $L^2(S)$. We have followed the argument from 
\cite{r190}, p.160.
\end{proof}

\begin{lemma}\label{s2_lm1}
Given $g\in L_2(S)$, let $w$ be the outgoing solution of the exterior
Dirichlet problem
$(\nabla^2 + k^2) w=0,\; \text{in}\ D^\prime$ with $w = g$ on $S$.
Then there exists a constant $C>0$, independent of $w$, such that
\be\label{s2_t9}
|||w|||\leq C\|g\|_{L^2(S)}\,,
\end{equation}
where $|||\cdot|||:=
||\cdot||_{H_{loc}^m(D')}+||\cdot||_{L^2(D';
(1+|x|)^{-\gamma})}$, $\gamma >1$, $m>0$ is an arbitrary integer, and
$H^m$ is the Sobolev space.

\end{lemma}

\begin{proof}
Let $G$ be the Dirichlet Green's function of the Laplacian in $D^\prime$:
\be\label{s2_w1}
\left(\nabla^2 + k^2 \right)G =-\delta (x-y) \hbox{\ in\ } D^\prime,
  \quad G=0 \hbox{\ on\ } S,
\end{equation}
\be\label{s2_w2}
\lim_{r \to \infty} \int_{|x| = r} \left|\frac{\partial G}{\partial |x|}
-  ik G \right|^2 ds = 0.
\end{equation}
Let $N$ be the unit normal to $S$ pointing into $D'$.
By Green's formula one has
\be\label{s2_w3}
w(x)= \int_S g(s) \frac{\partial G}{\partial N} (x,s) ds,\quad x\in
D'\,.
\end{equation}
The estimate for the $H_{loc}^m(D')$-norm part of \eqref{s2_t9}
follows from this representation and from the Cauchy inequality:
\[
|D^{(j)}w(x)|\leq ||g||_{L^2(S)}\left\| \frac{\partial
D_x^{(j)}G}{\partial N} (x,s)\right\|\leq c(x)||g||_{L^2(S)},
\]
where $c(x)\leq c(d)$ for all $x\in D'$ such that the distance $dist 
(x,S)\geq
d>0$.

For the $L^2$-weighted norm part of \eqref{s2_t9} let
 $R>0$ be such that $D\subset B_R=\{x\in\R^3\ :\ |x|<R\}$. Let $D'_R=B_R\setminus D$, and $S_R$
be the boundary of $B_R$.
The estimate
\be\label{s2_w4}
\left|\frac{\partial G}{\partial N} (x,s)\right| \leq \frac{c}{1 + |x|}, \quad |x| \geq R,
\end{equation}
and formula \eqref{s2_w3} imply 
\be\label{s2_w5} \|w\|_{L^2(S_R)}\leq c
\|g\|_{L^2(S)}, 
\end{equation} 
where here and in the sequel $c$ and $ C$ denote
various constants. Also, using the Cauchy inequality, formula  
\eqref{s2_w3}, inequality
\eqref{s2_w4} and the assumption $\gamma>1$, one gets 
\be\label{s2_w8}
\|w\|_{L^2(|x|>R;(1+\abs{x})^{- \gamma})}\leq c
\|g\|_{L^2(S)}\left\|\frac{1}{(1+|x|)^{\gamma+1}}\right\|_{L^2(|x|>R)}\leq
c \|g\|_{L^2(S)}. 
\end{equation}

To get the estimate for $\|w\|_{L^2(D'_R)}$ choose $R$ such that $k^2$ is not
a Dirichlet eigenvalue of $-\Delta$ in $D'_{R}$.
Then (\cite{lions}, p.189):
\be\label{s2_w10}
\|w\|_{H^m(D'_R)} \leq c [||(\Delta
+k^2)w ||_{\mathcal{H}^{m-2}(D'_R)}+
||w||_{H^{m-0.5}(S_R)} + ||w||_{H^{m-0.5}(S)}].
\end{equation}

The space $\mathcal{H}$ in the first term of the right-hand side in
\eqref{s2_w10} is different from the usual Sobolev space, but this term
is equal to zero anyway because $(\Delta +k^2)w =0$.

Let $m=0.5$ in \eqref{s2_w10}. Then
\be\label{s2_w11}
\|w\|_{H^{0.5}(D'_R)} \leq c [
||w||_{L^2(S_R)} + ||w||_{L^2(S)}].
\end{equation}

Since $w=g$ on $S$, then \eqref{s2_w5} and \eqref{s2_w11} imply
\be\label{s2_w12}
\|w\|_{L^2(D'_R)} \leq c ||g||_{L^2(S)}.
\end{equation}
\end{proof}

{\it Proof of Theorem \ref{s2_thm1}, continued}.

(ii)
Inequality \eqref{s2_t3} is the same as \eqref{s2_t1}, since $v=-u_0$ on $S$.
Estimate \eqref{s2_t4} follows from \eqref{s2_t3} and  Lemma \ref{s2_lm1}.

(iii)
Inequality \eqref{s2_t3} yields the convergence of $v_\epsilon$ to $v$ in the norm 
$\|\cdot\|_{L^2(S)}$. By  \eqref{s2_w5} $\|v_\epsilon-
v\|_{L^2(S_R)}\to 0$,
as $\epsilon \to 0$. On $S_R$ one has $v=\sum_{\l=0}^\infty 
A_\l(\alpha)\Psi_\l$ and $v_\epsilon=\sum_{\l=0}^{L(\epsilon)}c_\l 
\Psi_\l$. Multiply $v_\epsilon (R, \alpha')-v(R, \alpha')$ by
$\overline {Y_\l(\alpha')}$, integrate over $S^2$ and then let $\epsilon 
\to 0$. The result is (iii), and the proof of Theorem \ref{s2_thm1} is
completed.
\qed

The difference between RC and MRC is: (\ref{s2_t1}) does
not hold if one replaces $v_\epsilon$ by $\sum_{\l=0}^L
A_\l(\alpha)\Psi_\l$, and lets $L\to \infty$ (instead of
letting $\epsilon \to 0$). Indeed, the series
$\sum_{\l=0}^\infty A_\l(\alpha)\Psi_\l$ diverges at some
points of the boundary for many obstacles. Note also that
the coefficients in (\ref{s2_t2}) depend on $\epsilon$, so
(\ref{s2_t2}) is {\it not} a partial sum of a series.

For the Neumann boundary condition one minimizes
$$ \left\|\frac {\partial [u_0+\sum_{\l=0}^{L}c_\l\psi_\l]}{\partial
N}\right\|_{L^2(S)}$$
with respect to $c_\l$, and obtains essentially the same results.

According to Theorem \ref{s2_thm1} the computation of the outgoing
solution to (\ref{s1_e1})-(\ref{s1_e2}) is reduced to the
approximation of the boundary values in (\ref{s1_e1}) by the linear 
combinations of the functions $\Psi_\l$ 
restricted to the boundary $S$.
 A direct implementation of the above algorithm
is efficient for domains $D$ not very different from a circle, e.g.
for an
ellipse with a small eccentricity, but it fails for more complicated regions.
The numerical difficulties happen because the spherical Hankel functions
$h_l$ with large values of $l$ are bigger than $h_l$ with
small values of $l$ by many orders of magnitude. A finite precision of
numerical computations makes it necessary to
keep the values of $L$
not too high, e.g. $L\leq 20$. This restriction can be remedied
by the following modification of the above algorithm, see \cite{r437,r493}:

\begin{theorem}\label{s2_thm2}

Let $v:=u-u_0$, where $u$ is the solution to 
(\ref{s1_e1})-(\ref{s1_e2}).
Let $\ep>0$, and $L$ be a nonnegative integer. Suppose $U$ is an open
subset of $D$.

Then there exist a finite subset
$\{z_1,z_2,...,z_J\}\subset U$, and the coefficients 
$c_\l(\ep,z_j),\, 0\leq\l\leq L,\ 1\leq j\leq J=J(\ep),$ 
such that
the following inequalities \eqref{sw2_e61} and \eqref{s2_e622} hold:

\begin{itemize}
   \item [(i).]
\be\label{sw2_e61}
||u_0+v_\ep||_{L^2(S)} \leq \epsilon\,,
\end{equation}
where
\be\label{sw2_e621}
v_\ep(x):=\sum_{j=1}^J\sum_{\l=0}^L c_\l(\ep,z_j)\psi_\l(x,z_j),
\end{equation}
and
\begin{equation}\label{s2_hy}
 \psi_\l(x,z)= Y_\l(\alpha^\prime)h_\l(k|x-z|),\quad \alpha^\prime =
  \frac{x-z}{|x-z|},\quad z\in D,\quad x\in \R^3\setminus D.  
\end{equation}

   \item [(ii).]
\be\label{s2_e622}
\|v_\ep-v\|_{L^2(S)}\leq\ep
\end{equation}
and
\be\label{s2_e623}
|||v_\ep-v|||=O(\ep)\,,\quad \ep\ra 0\,,
\end{equation}
where
\[
|||\cdot|||=\|\cdot\|_{H^m_{loc}(D')}+\|\cdot\|_{L^2(D';(1+\abs{x})^{-
\gamma})}\,,
\]
$\gamma>1\,,\;m>0$ is an arbitrary integer, and $H^m$ is the 
Sobolev
space.

 \end{itemize}
\end{theorem}

\begin{proof}
(i) Note that in Theorem \ref{s2_thm1} we had $L=L(\ep)$, while now
we have $L$ fixed and $J=J(\ep)$. But the proof of Theorem 
\ref{s2_thm2} is similar to that of Theorem \ref{s2_thm1}.
 Let $\{z_j\}_{j=1}^\infty$ be a countable dense subset of $U$. To
establish \eqref{sw2_e61} it is sufficient to show that
\be\label{sw2_e63}
H:=\overline{span}\{\psi_\l(s,z_j)\ :\ 0\leq\l\leq L,\quad j=1,2,...\}=L^2(S)\,.
\end{equation}
Suppose that there exists $p\in L^2(S),\; p\not=0$ such that $p\perp H$ in $L^2(S)$.
Define the single-layer potential by
\be\label{sw2_e64}
W(y)=\int_S \frac{e^{ik|s-y|}}{ |s-y|} p(s)\ ds\,,\quad 
y\in\R^3\,.
\end{equation}
Then
\be\label{sw2_e65}
W(z_j)=\int_S \psi_0(s,z_j) p(s)\ ds=0
\end{equation}
for $j=1,2,...$.

The continuity of the single-layer potential in $\R^3$ implies that
$W(y)=0$ for all $y\in U$. The rest of the proof is the same as in
Theorem \ref{s2_thm1}.

\end{proof}

{\it Remark.} Functions $\{\Psi_l\}^\infty_{\l=0}$ are linearly independent on $S$. Indeed,
if some finite combination of these functions vanishes on
$S$, then it also vanishes in the exterior domain $D'$, since such a combination is an 
outgoing solution of the exterior Dirichlet problem with zero boundary conditions on $S$.
In particular, such a combination also vanishes on $S_R$. Since the spherical functions are orthogonal
on $S_R$, it implies that such a combination must be trivial.

See  Sections 6 and 7 for an extension of the MRC method to 
static problems,
and to scattering by periodic structures, respectively.

\section{Iterative MRC algorithms}\label{section_iter}

Let $z$ be a point in the interior of the obstacle $D$, and $x\in \R^3\setminus D$.
Recall that 
\begin{equation}\label{s3_hyy}
 \psi_\l(x,z)= Y_\l(\alpha^\prime)h_\l(k|x-z|), 
\end{equation}
where $h_\l (r)$ are the spherical
Hankel functions, normalized so that
$h_\l (r) \sim \frac{e^{ikr}}{r}$ as $r \to +\infty$.

{\bf Noniterative MRC.}

In this MRC implementation  one 
chooses a set of interior points $H=\{x_j\in D, j=1,2,...,J,\; J>0\}$ and minimizes

\be\label{s3_phi1}
 \Phi(\bc)=\|u_0(s)+\sum_{j=1}^J\sum_{\l=0}^L
c_{\l,j}\psi_\l(s,x_j)\|_{L^2(S)}, 
\end{equation}
 over
$\bc\in \C^N$, where $\bc=\{c_{\l,j}\}$. That is, the total
field $u(s)=u_0(s)+v(s)$ is desired to be as close to zero as
possible at the boundary $S$, to satisfy the required
condition for soft scattering. If the resulting residual
$r^{min}=\min\Phi$ is smaller than the prescribed tolerance
$\ep$, then the procedure is finished, and the sought
scattered field is
$$
v_{\epsilon}(x)=
\sum_{j=1}^J\sum_{\l=0}^{L}c_{\l,j}\psi_\l(x,x_j),\quad
x\in D'.
$$
If the residual $r^{min}>\ep$ then the method fails.
This approach, which can be called a {\bf Multi-point MRC},
 is justified by Theorem \ref{s2_thm2}. See \cite{r437,
r482, doicu} for details and results of numerical experiments.
The results show that the method is very efficient for domains $D$ of a
nearly spherical shape, i.e. without elongated parts. Clearly, the
only limitation in this method is the computer resources. The method
becomes impractical for large sets of interior points $H$.

To remedy this situation one can use iterative MRC implementations, of
which we describe the one based on a random choice of interior points, and
another one based on an optimal choice of such points.

{\bf Iterative  MRC with a random choice of points.}

Informally, the Random Multi-point MRC algorithm can be
described as follows.

Fix a $J>0$. Let $x_j, j=1,2,...,J$ be a batch of points randomly chosen  inside the
obstacle $D$. 

Let
$g(s)=u_0(s),\; s\in S$, and minimize the discrepancy
\be\label{s1_phi}
 \Phi(\bc)=\|g(s)+\sum_{j=1}^J\sum_{\l=0}^L
c_{\l,j}\psi_\l(s,x_j)\|_{L^2(S)} 
\end{equation}
 over
$\bc\in \C^N$, where $\bc=\{c_{\l,j}\}$. That is, the total
field $u(s)=g(s)+v(s)$ is desired to be as close to zero as
possible at the boundary $S$, to satisfy the required
condition for soft scattering. If the resulting residual
$r^{min}=\min\Phi$ is smaller than the prescribed tolerance
$\ep$, then the procedure is finished, and the sought
scattered field is
$$
v_{\epsilon}(x)=
\sum_{j=1}^J\sum_{\l=0}^{L}c_{\l,j}\psi_\l(x,x_j),\quad
x\in D',
$$

If, on the other hand, the residual $r^{min}>\ep$, then we
continue by trying to improve on the already obtained fit in
(\ref{s1_phi}). Adjust the field on the boundary by letting
$g(s):=g(s)+v_{\epsilon}(s),\; s\in S$. Create another batch
of $J$ points randomly chosen in the interior of $D$, and
minimize (\ref{s1_phi}) with this new $g(s)$. 
Continue with
the iterations until the required tolerance $\ep$ on the
boundary $S$ is attained. In each iteration accumulate
new interior points $x_j$ and the corresponding
best fit coefficients $c_{\l,j}$. After the desired tolerance is
reached, the sought scattered field $v_\ep$ is computed anywhere in $D'$.

Here is a precise description of the algorithm.

{\bf Random Multi-point MRC.}

For $x_j\in D$, and $\l\geq 0$ functions $\psi_\l(x,x_j)$
are defined as in (\ref{s3_hyy}).

\begin{enumerate} \item\label{init4} {\bf Initialization.}
Fix $\ep>0, \; L\geq 0,\;J>0,\; N_{max}>0$. Let $n=0$,
 and $g(s)=u_0(s), \;s\in S$.

\item\label{iter4} {\bf Iteration.}

\begin{enumerate}
 \item Let $n:=n+1$. Randomly choose $J$ points $x_j^{(n)}\in D,\; j=1,2,\dots, J$.

 \item
Minimize
\[
\Phi(\bc)=\|g(s)+\sum_{j=1}^J\sum_{\l=0}^L c_{\l,j}\psi_\l(s,x_j^{(n)})\|_{L^2(S)}
\]
over $\bc\in \C^N$, where $\bc=\{c_{\l,j}\}$.

Let the minimum of $\Phi$ be attained at
$\bc^{(n)}=\{c_{\l,j}^{(n)})\},\; j=1,2,\dots,J$, and 
the minimal value of $\Phi$ be $r^{min}$.

\end{enumerate}

\item {\bf Stopping criterion.}
\begin{enumerate}
 \item If $r^{min}\leq\epsilon$, then  stop.
 Compute the approximate scattered field anywhere in $D'$ by
\be\label{s1_vep}
v_{\epsilon}(x):=\sum_{k=1}^n \sum_{j=1}^J\sum_{\l=0}^{L}c_{\l,j}^{(k)}\psi_\l(x,x_j^{(k)}),\quad x\in D'.
\end{equation}

 \item If $r^{min}>\epsilon$, and $n\not=N_{max}$,
let
\[
g(s):=g(s)+\sum_{j=1}^J\sum_{\l=0}^{L}c_{\l,j}^{(n)}\psi_\l(s,x_j^{(n)}),\quad x\in S
\]
and repeat the iterative
step (\ref{iter4}).
\item If $r^{min}>\ep$, and $n=N_{max}$, then the procedure failed.

\end{enumerate}
\end{enumerate}

Numerical experiments based on this method are presented in the next
section. The method is relatively slow, and it can be improved by
choosing the interior points in some optimal way.

{\bf Iterative MRC with an optimal choice of points.}

In this case the interior points $z_1, z_2,...$ in $D$ are chosen one at a time, and 
their placement is not random. Rather, the discrepancy $\Psi$ is
minimized not only with respect to the coefficients $\bf c$, but also
with respect to the position of these points $z_j$.

 Let
$g_1(s)=u_0(s)=u_0(s, \alpha),\; s\in S$.

Minimize
 \be\label{so1_phi}
\Phi(z_1,\bc(z_1)):=\min_{z\in D}\min_{\bc\in
\C^N}\|g_1(s)+\sum_{\l=0}^L c_{\l}\psi_\l(s,z)\|_{L^2(S)},
\end{equation}
where $\bc=\{c_{\l}\}=\{c_{\l m}\}_{0\leq \l \leq L, -\l\leq m\leq \l}$, 
$L\geq 0$ is a fixed integer, and
$\sum_{\l=0}^L:=\sum_{\l=0}^L\sum_{m=-\l}^{\l}$.
 Let
 \be\label{so1_v}
v_1(x)=\sum_{\l=0}^L c_{\l}(z_1)\psi_\l(x,z_1), \quad 
c_{\l}(z_1)=c_{\l}(z_1, \alpha).
\end{equation}
 The
requirement \eqref{so1_phi}
 means that the total field $u(s)=g_1(s)+v_1(s)$ has to be as close
to zero as possible on the boundary $S$, so that it approximates
best the Dirichlet boundary condition in \eqref{s1_e1}. This is
achieved by varying the interior point $z\in D$ and choosing the
coefficients $\bc(z)\in\C^N$ giving $g_1+v_1$ the best fit to zero on the boundary
$S$.
 Let the minimum in \eqref{so1_phi} be attained at $z_1\in D$.
If the resulting value of the residual
$r^{min}=\Phi(z_1,\bc(z_1))$ is smaller than the prescribed
tolerance $\ep$, than the procedure is finished. The sought
approximate scattered field is $v_1(x),\;x\in D'$ (see Theorem
\ref{s2_thm2}), and the approximate scattering amplitude is
\be\label{so1_a} A_1(\a',\a)= e^{-ik\a'\cdot z_1}
\sum_{\l=0}^{L}c_{\l}(z_1)Y_\l(\a') \,.
\end{equation}
Note that $c_{\l}(z_1)=c_{\l}(z_1,\a)$.

The expression for $A_1(\a',\a)$ in \eqref{so1_a} is obtained from
\eqref{so1_v} by letting $|x|\ra\infty$ in $x=\a'|x|$, because of
our normalization
 \be\label{so1_h}
 h_\l(k|x|)=\frac {e^{ik|x|}}{|x|}
\left\{1+O\left(\frac 1{|x|}\right)\right\},\quad |x|\ra\infty,
\end{equation}
and $|x-z|=|x|-\a'\cdot z+O(1/|x|)$ as $|x|\ra\infty$.

If, on the other hand, the residual $r^{min}>\ep$, then we
continue by trying to improve on the already obtained fit in
(\ref{so1_phi}) as follows. Adjust the field on the boundary by
letting $g_2(s)=g_1(s)+v_1(s),\; s\in S$, and do the minimization
\eqref{s1_phi} with $g_2(s)$ instead of $g_1(s)$, etc. Continue
with the iterations until the required tolerance $\ep$ on the
boundary $S$ is attained. At the same time keep track of the
changing approximate scattered field $v_n(x)$, and the scattering
amplitude $A_n(\a',\a)$. In this construction $g_{n+1}=u_0+v_n$ on
$S$. The goal of \eqref{s1_phi} is to obtain $g_n\ra 0$  in $L^2(S)$
as $n\ra\infty$, yielding
$u_0+v_n\ra 0$ in $L^2(S)$. According to Theorem \ref{s2_thm2}, this gives 
an
approximate scattered solution $v_n$ on $D'$ to \eqref{s1_e1}-\eqref{s1_e2}.

Here is a precise description of the algorithm.
\medskip

{\bf MRC method with optimal choice of sources.}

\begin{enumerate}
 \item\label{init5}
{\bf Initialization.} Fix $\ep>0, \; L\geq 0,\; N_{max}>0$. Let
$n=0,\;v_0(x)=0,\;A_0(\a',\a)=0$, and $g_1(s)=u_0(s), \;s\in S$.

\item\label{iter5} {\bf Iteration.}

\begin{enumerate}
\item Increase the value of $n$ by $1$.

 \item
Minimize
\[
\Phi(z_n,\bc(z_n)):=\min_{z\in D}\min_{\bc\in
\C^N}\|g_n(s)+\sum_{\l=0}^L c_{\l}\psi_\l(s,z)\|_{L^2(S)},
\]

with the minimal value attained at $z_n\in D,\;\bc(z_n)\in\C^N$.

 \item Let
\[
v_n(x)=v_{n-1}(x)+\sum_{\l=0}^{L}c_{\l}(z_n)\psi_\l(x,z_n),\quad
x\in D',
\]
\[
A_n(\a',\a)=A_{n-1}(\a',\a)+e^{-ik\a'\cdot z_n}
\sum_{\l=0}^{L}c_{\l}(z_n)Y_\l(\a'),
\]
and
\[
g_{n+1}(s)=g_n(s)+\sum_{\l=0}^{L}c_{\l}(z_n)\psi_\l(s,z_n),\quad
s\in S,
\]
that is $g_{n+1}(s)=u_0(s)+v_n(s),\;s\in S$.
 \item Let
\[
r^{min}:=\Phi(z_n,\bc(z_n)).
\]
\end{enumerate}

\item {\bf Stopping criterion.}
\begin{enumerate}
 \item If $r^{min}\leq\epsilon$, then  stop; $v_n(x)$ is the approximate
scattered field, and $A_n(\a',\a)$ is the approximate scattering
amplitude.
 \item If $r^{min}>\epsilon$, and $n<N_{max}$, then repeat the iterative
step (\ref{iter5}).
\item If $r^{min}>\ep$, and $n=N_{max}$, then the procedure failed.

\end{enumerate}
\end{enumerate}

\section{Numerical Experiments for Random Multi-point MRC}

In this section we describe numerical results obtained by
the Random Multi-point MRC method for 2D and 3D obstacles.
We also compare the 2D results
 to the ones obtained by the Multiple-point MRC described above, and introduced in \cite{r437}. 
The results of \cite{r437} show a favorable comparison of
the Multi-point MRC method with the Boundary Integral Equation Method. Further improvements are attained
with the Random Multi-point MRC method, for which the Multi-point MRC is just the first iteration.

Note that in a 2D case instead of (\ref{s3_hyy}) one has
\[
\psi_l(x,x_j)=H_l^{(1)}(k\abs{x-x_j})e^{il\theta_j},
\]
where $(x-x_j)/\abs{x-x_j}=e^{i\theta_j}$.

For a  numerical implementation choose
$M$ nodes $\{t_m\}$ on the surface $S$ of the obstacle $D$. 
After the interior points 
$x_j,\; j=1,2,...,J$ are chosen, form   $N$ vectors 
\[
{\bf a}^{(n)}=\{\psi_l(t_m,x_j)\}_{m=1}^M,
\]
$n=1,2,\dots,N$ of length $M$. Note that $N=(2L+1)J$
for a 2D case, and $N=(L+1)^2J$ for a 3D case.
It is convenient to normalize the norm in $\R^M$ by
\[
\|{\bf b}\|^2=\frac 1M \sum_{m=1}^M|b_m|^2,\quad {\bf
b}=(b_1,b_2,...,b_M).
\]
Then $\|u_0\|=1$.

Now let ${\bf b}=\{g(t_m)\}_{m=1}^M$, in
the Random Multi-point MRC (see Section \ref{section_iter}), and minimize
\be\label{s3_minm}
\Phi({\bf c})=\|{\bf b}+A{\bf c}\|,
\end{equation}
for ${\bf c}\in \C^N$, where $A$ is the matrix containing vectors ${\bf a}^{(n)},\;
n=1,2,\dots,N$ as its columns.

The Singular Value Decomposition (SVD) method (see e.g.
\cite{numrec}) is used to minimize (\ref{s3_minm}). 
According to the SVD, the matrix $A$ is represented as
\[A=UWV^H\,,
\]
where the $M\times N$ matrix $U$ has orthonormal columns ${\bf u}^{(n)}, \; n=1, \ldots, N$,
the square $N\times N$ matrix $V$ has orthonormal columns ${\bf v}^{(n)}$, $n=1, \ldots, N$,
and the diagonal square
$N\times N$ matrix $W=(w_n)^N_{n=1}$ is composed of the (nonnegative) singular values of $A$.

Let $P\subset\{1,2,\dots, N\}$ be defined by
\[
P=\{n\ :\ w_n\geq w_{min}\}\,
\]
for some positive constant $w_{min}$.

 Compute the normalized residual
\[
r^{min}=\frac 1{\sqrt{M}} \sqrt{ \|{\bf b}\|^2-\sum_{n\in P} \abs{<{\bf u}^{(n)},{\bf
b}>}^2}\,.
\]

The minimizer ${\bf c}$ is given by
\[
{\bf c}=\sum_{n\in P} \frac 1{w_n}<{\bf u}^{(n)},{\bf b}>{\bf
v}^{(n)}\,.
\]

 Small singular values $w_n<w_{min}$ of the
matrix $A$ are used to identify and delete linearly dependent or almost
linearly dependent combinations of vectors ${\bf a}^{(n)}$. This spectral
cut-off makes the minimization process stable, see details in
\cite{r437}. 

 Let $r^{min}$ be the residual, i.e. the minimal value 
of
$\Phi({\bf c})$ attained after $N_{max}$ iterations of the
Random Multi-point MRC method (or when it is stopped). For a
comparison, let $r^{min}_{old}$ be the residual obtained in
\cite{r437} by the Multi-point MRC.

We have conducted 2D numerical experiments for four obstacles: two 
ellipses of
different eccentricity, a kite, and a triangle. 
The M=720 nodes $t_m$ were uniformly distributed on the interval
$[0,2\pi]$, used to parametrize the boundary $S$.
Each case was tested for
wave numbers $k=1.0$ and $ k=5.0$. Each obstacle was subjected to incident
waves corresponding to $\a=(1.0,0.0)$ and $\a=(0.0,1.0)$. 

The results for the Random Multi-point MRC with $J=1$ are shown in Table 1, in the
last column $r^{min}$. In every experiment the target
residual $\ep=0.0001$ was obtained in under 6000 iterations,
 in about 2 minutes run time on a 2.8 MHz PC.

In \cite{r437}, we have conducted numerical experiments for the same four 
2D
obstacles by a Multi-point MRC, as described in the beginning of this
section. The interior points $x_j$ were chosen differently in
each experiment. Their choice is indicated in the description of each 2D
experiment.
The column $J$ shows the number of these interior
points. Values $L=5$ and $M=720$ were used in all the
experiments. These results are
shown in Table 1, column  $r^{min}_{old}$.

 Thus, the Random Multi-point MRC method achieved a significant
improvement over the Multi-point MRC. 

\begin{table} 
\caption{Normalized residuals attained in the
numerical experiments for 2D obstacles by Random Multi-point MRC,
 $\|{\bf u_0}\|=1$.}

\begin{center}
\begin{tabular}{c  r  r  c  r   r}

\hline
Experiment & $J$ & $k$ & $\a$ & $r^{min}_{old}$ & $r^{min}$ \\

\hline
I & 4 & 1.0 & $(1.0,0.0)$ & 0.000201  & 0.0001\\
Ellipse  & 4 & 1.0 & $(0.0,1.0)$ & 0.000357  & 0.0001\\
  & 4 & 5.0 & $(1.0,0.0)$ & 0.001309  & 0.0001\\
  & 4 & 5.0 & $(0.0,1.0)$ & 0.007228  & 0.0001\\
\hline
II & 16 & 1.0 & $(1.0,0.0)$ & 0.003555  & 0.0001\\
Kite  & 16 & 1.0 & $(0.0,1.0)$ & 0.002169  & 0.0001\\
  & 16 & 5.0 & $(1.0,0.0)$ & 0.009673  & 0.0001\\
  & 16 & 5.0 & $(0.0,1.0)$ & 0.007291  & 0.0001\\
\hline
III & 16 & 1.0 & $(1.0,0.0)$ & 0.008281  & 0.0001\\
Triangle    & 16 & 1.0 & $(0.0,1.0)$ & 0.007523  & 0.0001\\
    & 16 & 5.0 & $(1.0,0.0)$ & 0.021571  & 0.0001\\
    & 16 & 5.0 & $(0.0,1.0)$ & 0.024360  & 0.0001\\
\hline
IV & 32 & 1.0 & $(1.0,0.0)$ & 0.006610  & 0.0001\\
Ellipse   & 32 & 1.0 & $(0.0,1.0)$ & 0.006785  & 0.0001\\
   & 32 & 5.0 & $(1.0,0.0)$ & 0.034027  & 0.0001\\
   & 32 & 5.0 & $(0.0,1.0)$ & 0.040129  & 0.0001\\
\hline
\end{tabular}
\end{center}              
\end{table}

{\bf Experiment 2D-I.} The boundary $S$ is an ellipse described by
\be
{\bf r}(t)=(2.0\cos t,\ \sin t),\quad 0\leq t<2\pi\,.
\end{equation}
The Multi-point MRC used $J=4$ interior points 
$x_j=0.7{\bf r}(\frac{\pi(j-1)}2),\; j=1,\dots,4$.
The run time was 2 seconds.

{\bf Experiment 2D-II.} The kite-shaped boundary $S$ (see \cite{coltonkress}, Section 3.5) is  
described by
\be
{\bf r}(t)=(-0.65+\cos t+0.65\cos 2t,\ 1.5\sin t),\quad 0\leq t<2\pi\,.
\end{equation}    
The Multi-point MRC used $J=16$ interior points $x_j=0.9{\bf r}(\frac{\pi(j-1)}8),\; j=1,\dots,16$.
The run time  was 33 seconds.


{\bf Experiment 2D-III.} The boundary $S$ is the triangle
with vertices at $(-1.0,0.0)$ and $(1.0,\pm 1.0)$.
The Multi-point MRC used the interior points $x_j=0.9{\bf r}(\frac{\pi(j-1)}8)$,
$ j=1,\dots,16$.
The run time  was about 30 seconds.

{\bf Experiment 2D-IV.} The boundary $S$ is an ellipse described by
\be
{\bf r}(t)=(0.1\cos t,\ \sin t),\quad 0\leq t<2\pi\,.
\end{equation}
The Multi-point MRC used $J=32$ interior points 
$x_j=0.95{\bf r}(\frac{\pi(j-1)}{16}),\; j=1,\dots,32$.
The run time was about 140 seconds.

The 3D numerical experiments were conducted for 3 obstacles: a sphere, a cube, and an ellipsoid.
We used the Random Multi-point MRC with $L=0,\; w_{min}=10^{-12}$, and $J=80$.
 The number $M$ of the points on the boundary $S$ is
indicated in the description of the obstacles. The scattered
field for each obstacle was computed for two incoming
directions $\a_i=(\theta,\phi),\;i=1,2$, where $\phi$ was
the polar angle. The first unit vector $\a_1$ is denoted by
(1) in Table 2, $\a_1=(0.0,\pi/2)$. The second one is
denoted by (2), $\a_2=(\pi/2,\pi/4)$. A typical number of
iterations $N_{iter}$ and the run time on a 2.8 MHz PC are
also shown in Table 2. For example, in experiment I with
$k=5.0$ it took about 700 iterations of the Random
Multi-point MRC method to achieve the target residual
$r^{min}=0.001$ in 7 minutes.

{\bf Experiment 3D-I.} The boundary $S$ is the sphere of radius $1$,
with $M=450$.

{\bf Experiment 3D-II.} The boundary $S$ is the surface of the cube
$[-1,1]^3$ with $M=1350$.

{\bf Experiment 3D-III.} The boundary $S$ is the surface of the
ellipsoid $x^2/16+y^2+z^2=1$ with $M=450$.

\begin{table}
\caption{Normalized residuals attained in the numerical experiments for 3D obstacles
 by Random Multi-point MRC,
$\|{\bf u_0}\|=1$.}

\begin{center}

\begin{tabular}{c  r  r  l  r   r}

\hline
Experiment & $k$ & $\a_i$ &  $r^{min}$ & $N_{iter}$ & run time \\

\hline
I & 1.0 &  & $0.0002$ & 1   & 1 sec\\
Sphere  & 5.0 &  & $0.001$ & 700  & 7 min\\

\hline
II & 1.0 & (1) & $0.001$ & 800   & 16 min\\
Cube   & 1.0 & (2) & $0.001$ & 200  & 4 min\\
   & 5.0 & (1) & $0.0035$ & 2000   & 40 min\\
   & 5.0 & (2) & $0.002$ & 2000  & 40 min\\

\hline
III & 1.0 & (1) & $0.001$ & 3600   & 37 min\\
Ellipsoid    & 1.0 & (2) & $0.001$ & 3000   & 31 min\\
    & 5.0 & (1) & $0.0026$ & 5000   & 53 min\\
    & 5.0 & (2) & $0.001$ & 5000   & 53 min\\
\hline
\end{tabular}
\end{center}

\end{table} 
In the last experiment the run time could be reduced by taking a smaller value
for $J$. For example, the choice of $J=8$ reduced the running time to
about 6-10 minutes.

Numerical experiments show that the minimization results depend on the
choice of such parameters as $J,\; w_{min}$, and $L$.

\section{Numerical Experiments for Optimal Choice MRC}

In this section we describe numerical results obtained by the  MRC
method with the optimal choice of sources for 2D and 3D obstacles.
The notations are kept the same as in the previous section for the Random Multi-point MRC.
As there, one has to minimize

\be\label{s4_minm} \Phi({z,\bf c})=\|{\bf b}+A{\bf c}\|,
\end{equation}
in every iterative step, but, in addition,  the
residual is minimized with respect to the interior point $z\in D$.

There is a variety of methods to minimize $\Phi(z,\bc(z))$, since after
the minimization in the coefficients $\bc(z)$ by the SVD
it is just a 2D or 3D minimization in the region $D$. Our choice was
the Powell's method which imitates the conjugate gradients approach,
but does not require analytical expressions for the gradient. The
Brent  method was used for a line minimization, see \cite{numrec,r493} for details.
The Powell's algorithm is also described below in Section
\ref{section_sim}.

\begin{table}[tb]
\caption{Normalized residuals attained in the
numerical experiments for 2D obstacles by the Optimal Choice MRC,
 $\|{\bf u_0}\|=1$.}

\begin{center}

\begin{tabular}{c  c  c  r  c   c}

\hline
Experiment &  $k$ & $\a$ & $N_{iter}$  & $r^{min}$ & \tiny{ (MRC-BIEM)/BIEM} \\

\hline
I &  1.0 & $(1.0,0.0)$ & 20 & 0.0010  & 0.0001\\
Ellipse   & 1.0 & $(0.0,1.0)$& 20 & 0.0018  & 0.0001\\
   & 5.0 & $(1.0,0.0)$& 53 & 0.0010  & 0.0001\\
   & 5.0 & $(0.0,1.0)$& 45 & 0.0020  & 0.0001\\
\hline
II  & 1.0 & $(1.0,0.0)$& 53 & 0.0020  & 0.0001\\
Kite     & 1.0 & $(0.0,1.0)$& 32 & 0.0020  & 0.0001\\
     & 5.0 & $(1.0,0.0)$& 75 & 0.0020  & 0.0003\\
     & 5.0 & $(0.0,1.0)$& 68 & 0.0020  & 0.0001\\
\hline
III  & 1.0 & $(1.0,0.0)$& 55 & 0.0020  & \\
Triangle    & 1.0 & $(0.0,1.0)$& 48 & 0.0017  & \\
    & 5.0 & $(1.0,0.0)$& 72 & 0.0019  & \\
    & 5.0 & $(0.0,1.0)$& 80 & 0.0020  & \\
\hline
IV  & 1.0 & $(1.0,0.0)$& 100 & 0.0041  & 0.0008\\
Ellipse   & 1.0 & $(0.0,1.0)$& 100 & 0.0027  & 0.0000\\
   & 5.0 & $(1.0,0.0)$& 100 & 0.0058  & 0.0004\\
   & 5.0 & $(0.0,1.0)$& 100 & 0.0037  & 0.0012\\
\hline
V &  1.0 & $(1.0,0.0)$ & 1 & 0.0000  & 0.0001\\
Circle   & 5.0 & $(1.0,0.0)$& 21 & 0.0020  & 0.0001\\
\hline
\end{tabular}
\end{center}

\end{table}

In addition to the four obstacles considered for the Random Multi-point
MRC, the circle $|r|=1$ was tested to check if the Optimal point
MRC was able to find the scattered field just after one iteration, since, 
in this case, the optimal point was in the origin. The result is in Table
3, experiment number V.
 The column
$N_{iter}$ shows the number of iterations (number of source
points) at the end of the iterative process. The process was
stopped after the algorithm reached the sought tolerance
$\ep=0.002$, or $N_{max}=100$.
Values $L=5$ and $M=720$ were used
in all 2D experiments.

The last column $(MRC-BIEM)/BIEM$ shows the discrepancy in the
scattering amplitude computed by the MRC and BIEM methods. The
values shown are the $L_2$ norms of the difference of the
scattering amplitude obtained by MRC and BIEM, over the $L_2$ norm
of the scattering amplitude obtained by BIEM. We followed
\cite{coltonkress} for the BIEM implementation using 64 points on
the boundary $S$ in every 2D experiment. No comparison is provided
for a triangular obstacle, since it requires a complete rewriting
of the BIEM code to accommodate the corner points. No such
rewriting is required for the MRC method. Table III shows that for the value
of tolerance $\ep=0.002$ the computed scattering amplitude is in an excellent agreement
with the scattering amplitude computed using BIEM. Results for 3D obstacles are 
provided in Table 4.

Concerning the efficiency of the methods: for simple 
geometries
the Multi-point MRC (see \cite{r437}) is the fastest, provided
that the required accuracy can be achieved by a relatively small
number $J$ of the interior points (sources) used simultaneously.
 This assures the
resulting matrices being of a manageable size. Otherwise, one has
to use Random, or Optimal choice MRC, which
take a significantly longer time to run, but can accomplish the
solution of scattering problems untractable by single step
methods, such as the Multi-point MRC or BIEM. While the precision of the
Random-point MRC was higher in the conducted experiments, the optimally placed
MRC method achieves an order of magnitude improvement in run time over
it.

\begin{table}[tb]
\caption{Normalized residuals attained in the numerical experiments
 for 3D obstacles by the Optimal Choice MRC,
$\|{\bf u_0}\|=1$.}
\begin{center}

\begin{tabular}{c  c  c  r  c  }

\hline
Experiment & $k$ & $\a_i$ &  $N_{iter}$ &$r^{min}$  \\

\hline
I & 1.0 & & 1 & $0.0000$   \\
Sphere  & 5.0 & & 43 & $0.0019$  \\

\hline
II & 1.0 & (1) & 12 & $0.0019$    \\
Cube   & 1.0 & (2)& 7 & $0.0019$   \\
   & 5.0 & (1) & 70 & $0.0019$    \\
   & 5.0 & (2) & 35 & $0.0020$   \\
\hline
III & 1.0 & (1) & 12 & $0.0016$   \\
Ellipsoid    & 1.0 & (2)  &  35  & $0.0020$  \\
    & 5.0 & (1) & 55 & $0.0020$   \\
    & 5.0 & (2) & 67 & $0.0020$    \\
\hline
\end{tabular}
\end{center}

\end{table}

\section{MRC for static problems}

In this Section we follow \cite{r481} and \cite{r476}, 
Chapter 12.
Consider a bounded domain $D \subset \R^3$ with a Lipschitz
boundary $S$, $D\subset B_R := \{x : |x| \leq R\}$. Denote
the exterior
 domain by $D^\prime = \R^3 \backslash D$. 
Let $ S^2$ denote the unit sphere in $\R^3$. 
Consider the problem:
\be\label{s6_11}
\nabla^2 v=0 \hbox{\ in\ } D^\prime, \quad
  v = f \hbox{\ on\ } S,
\end{equation}
\be\label{s6_12}
v:=O\left(\frac 1{r}\right),
 \quad r:=|x| \to \infty. 
\end{equation}

 Denote
by $Y_\l (\alpha),\;\alpha\in S^2 $  the orthonormal spherical harmonics,
$Y_\l = Y_{\l m},\\ -\l \leq m \leq \l$, 
and let harmonic functions $H_\l(x)$ be defined by
\[
H_\l(x) :=\frac 
{Y_\l(\alpha)}{r^{\l 
+1}},\quad\l\geq 0,\quad \alpha:=\frac x r\in S^2.
\]

In the region $r> R$ the solution to \eqref{s6_11}-\eqref{s6_12} is:
\be\label{s6_13}
v(x) = 
\sum^\infty_{\l =0} c_\l H_\l(x), \quad
\quad r > R.
\end{equation}
The summation in \eqref{s6_13} and below includes summation with respect 
to $m$, 
$-\l \leq m \leq \l$,
and $c_\l=c_{\l, m} $ are some coefficients determined by $f$.

The series \eqref{s6_13}  in general does not converge up to the 
boundary $S$.
Our aim is to give a formulation of an analog of the Modified Rayleigh 
Conjecture (MRC) from \cite{r430},
which can be used in numerical solution of the 
boundary-value  problems. 
The authors hope that the MRC method
for static problems can be used as a basis for an efficient
numerical algorithm for solving boundary-value problems for Laplace
equations in domains with complicated boundaries. In above sections such
algorithms were developed on the basis of MRC  for solving
boundary-value problems for the Helmholtz equation.
Although the boundary integral equation methods and finite elements 
methods are widely and successfully used for solving these problems, the 
method, based on 
MRC, proved to be competitive and often superior to the currently 
used methods. 

We discuss the Dirichlet condition but 
a similar argument is applicable to the Neumann and 
Robin boundary 
conditions. Boundary-value problems and scattering problems in 
rough 
domains were studied in \cite{rammb1} and in \cite{r476}, 
Chapter 9. 
  
Let us present the basic results on which the MRC method is based.

Fix $\epsilon >0$, an arbitrary small number.

{\bf Lemma 6.1.} {\it There exist $L=L(\epsilon)$ and 
$c_\l=c_\l(\epsilon)$
such that }
\be\label{s6_14}
 ||\sum_{\l=0}^{L(\epsilon)}c_\l(\epsilon)H_\l -f||_{L^2(S)} \leq 
\epsilon.
\end{equation}

If \eqref{s6_14} and the boundary condition \eqref{s6_11} hold, then
\be\label{s6_15}
||v_{\epsilon}-v||_{L^2(S)}\leq \epsilon,  \quad 
v_{\epsilon}:=\sum_{\l=0}^{L(\epsilon)}c_\l(\epsilon)H_\l.
\end{equation}

{\bf Lemma 6.2.} {\it If \eqref{s6_14} holds then
\be\label{s6_16}
||v_{\epsilon}-v||=O(\epsilon) \quad \epsilon \to 0,
\end{equation}
where $||\cdot||:= ||\cdot||_{H_{loc}^m(D')}+||\cdot||_{L^2(D'; 
(1+|x|)^{-\gamma})}$, $\gamma >1$, $m>0$ is an arbitrary integer,
and $H^m$ is the Sobolev space.} 

In particular, \eqref{s6_16} implies
\be\label{s6_17}
 ||v_{\epsilon}-v||_{L^2(S_R)}=O(\epsilon) \quad  \epsilon \to 0. 
\end{equation}

Let us formulate an analog of the Modified Rayleigh Conjecture (MRC):  

{\bf Theorem 6.1 (MRC):} {\it For an arbitrary small $\epsilon>0$ there 
exist
$L(\epsilon)$ and $c_\l(\epsilon), 0\leq \l \leq L(\epsilon)$,
such that \eqref{s6_14} and \eqref{s6_14}  hold.}

Theorem 6.1 follows from  Lemmas 6.1 and 6.2.

For the Neumann boundary condition one minimizes
$ ||\frac {\partial [\sum_{\l=0}^{L}c_\l\psi_\l]}{\partial 
N}-f||_{L^2(S)}$
with respect to $c_\l$. Analogs of Lemmas 6.1-6.2 are valid and their 
proofs are essentially the same.

If the boundary data $f\in C(S)$, then one can use $C(S)-$ norm 
in \eqref{s6_14}-\eqref{s6_17}, and an analog of Theorem 6.1 then 
follows immediately from the maximum principle.

To solve problem \eqref{s6_11}-\eqref{s6_12}
using MRC, fix a small $\epsilon >0$
and find $L(\epsilon)$ and $c_\l(\epsilon)$ such that 
\eqref{s6_14} holds. This is possible by Lemma 6.1 and can be done numerically
by minimizing $||\sum_0^Lc_\l 
H_\l -f||_{L^2(S)}:=\phi 
(c_1,.....,c_L)$. If the minimum of $\phi$ is larger than $\epsilon$, 
then increase $L$ and repeat the minimization. Lemma 6.1 guarantees the
existence of such $L$ and $c_\l$ that the minimum is less than $\epsilon$.
Choose the smallest $L$ for which this happens and
define $v_\epsilon:=\sum^L_{\l = 0} c_\l H_\l$. Then, by Lemma 6.2,
$v_\epsilon$ is the approximate solution to problem \eqref{s6_11}-\eqref{s6_12}
with the accuracy $O(\epsilon)$ in the norm $||\cdot||$.

{\bf Proof of Lemma 6.1.}   We start with the claim:

{\bf Claim}: 
  {\it The restrictions of harmonic functions $H_\l$ on $S$
form a total set in $L^2(S)$.} 

Lemma 6.1 follows from this claim.
Let us prove the claim. Assume the contrary. Then there is
a function $g\neq 0$ such that $\int_S g(s)h_\l(s) ds=0\,\, \forall \l\geq 
0.$ This implies $V(x):=\int_S g(s)|x-s|^{-1}ds=0\,\, \forall x\in D'$. 
Thus $V=0$ on $S$, and since $\Delta V=0$ in $D$, one concludes that $V=0$
in $D$. Thus $g=0$ by the jump formula for the normal derivatives of
the simple layer potential $V$. This contradiction  proves
the claim. Lemma 6.1 is proved. 
$\Box$

{\bf Proof of Lemma 6.2.} 
By Green's formula one has
\be\label{s6_31}
 w_\epsilon (x) = \int_S w_\epsilon (s) G_N (x,s) ds, \quad
  \| w_\epsilon  \|_{L^2(S)} < \epsilon,\quad  w_\epsilon:=v_\epsilon-v.
\end{equation}
Here $N$ is the unit normal to $S$, pointing into $D'$, and $G$ is the 
Dirichlet Green's function of the Laplacian in $D^\prime$:
\be\label{s6_32}
\nabla^2 G =-\delta (x-y) \hbox{\ in\ } D^\prime,
  \quad G=0 \hbox{\ on\ } S,
\end{equation}
\be\label{s6_33}
 G = O\left(\frac 1 r\right), \quad r\to \infty.
\end{equation}
 From \eqref{s6_31} one gets \eqref{s6_17} and \eqref{s6_16} 
with 
$H^m_{loc}(D')-$norm immediately by the 
Cauchy inequality. Estimate \eqref{s6_16} in the region 
$B_{R}^\prime:=\R^3\setminus B_R$ follows
from the estimate
\be\label{s6_34}
\left|G_N (x,s)\right| \leq \frac{c}{1 + |x|}, \quad |x| 
\geq R.
\end{equation}
 In the region $B_R\backslash D$ estimate \eqref{s6_16} follows from local elliptic 
estimates for $w_\epsilon:=v_\epsilon -v$, which 
imply that
\be\label{s6_35}
\|w_\epsilon \|_{L^2(B_R\backslash D)} \leq c\epsilon. 
\end{equation}
Let us recall the elliptic estimate we have used.
Let  $D'_{R}:=B_R\backslash D$ and $S_R$ be the boundary of $B_R$.
 Recall the elliptic estimate  for the solution to
homogeneous Laplace equation in $D'_{R}$ ( see \cite{lions}, p.189):
\be\label{s6_36}
\|w_\epsilon \|_{H^{0.5}(D'_{R})} \leq c [ 
||w_\epsilon||_{L^2(S_R)} + ||w_\epsilon||_{L^2(S)}].
\end{equation}
The estimates 
$||w_\epsilon||_{L^2(S_R)}=O(\epsilon)$, 
 $||w_\epsilon||_{L^2(S)}=O(\epsilon)$, and \eqref{s6_36} yield
\eqref{s6_16}. Lemma 6.2 is proved. \qed

\section{MRC for scattering by periodic structures}

Determination of fields scattered by periodic structures is of a great
importance in modern diffractive optics, and there is a vast literature
on both the direct and inverse problems of
this type, see, for example \cite{petit}. Still, an efficient computation of such fields presents
certain difficulties.
In the next Sections we present some theoretical background, a modification
of the MRC method, and numerical results for such a scattering, see
\cite{r461}.

For simplicity we consider a 2-D setting, but our arguments can be as
easily  applied to $n$-dimensional problems, $n\geq 2$.
Let $f : \rc\ra\rc,\ f(x+L)=f(x)$
be an $L$-periodic Lipschitz continuous function, and let $D$ be the
domain
\[
D=\{(x,y)\ :\ y\geq f(x),\ x\in \rc\}.
\]

Without loss of generality we assume that $f\geq 0$. If it is not, one can
choose the origin so that this assumption is satisfied, because
$M:=\sup_{0\leq x \leq L} |f(x)|<\infty$.

Let $\bx=(x,y)$ and $u(\bx)$ be the total field satisfying
\begin{equation}\label{s1_1}
(\Delta +k^2)u=0,\quad \bx\in D, \quad k=const>0
\end{equation}
\begin{equation}\label{s1_2}
u=0 \quad \text{on} \quad S:\,=\D D,
\end{equation}
\begin{equation}\label{s1_3}
u=u_0+v,\quad u_0:\,=e^{ik\a\cdot\bx},
\end{equation}
where the unit vector $\a=(\cos\theta,-\sin\theta),\ 0<\theta<\pi/2$,
and $v(\bx)$ is the scattered field, whose asymptotic behavior as
$y\ra\infty$ will be specified below, and
\begin{equation}\label{s1_4}
u(x+L,y)=\nu u(x,y),\quad u_x(x+L,y)=\nu u_x(x,y)\;\text{in}\; D,\quad
\nu:\,=e^{ikL\cos\theta}\,.
\end{equation}

Conditions (\ref{s1_4}) are the $qp$ ({\bf quasiperiodicity}) conditions.
To find the proper radiation condition for the scattered field $v(\bx)$
consider the spectral problem
\begin{equation}\label{s1_5}
\p''+\l^2\p=0, \quad 0<x<L,
\end{equation}
\begin{equation}\label{s1_6}
\p(L)=\nu\p(0), \quad \p'(L)=\nu\p'(0)
\end{equation}
arising from the separation of variables in (\ref{s1_1})-(\ref{s1_4}).
This problem has a discrete spectrum, and its eigenfunctions form a
basis in $L^2(0,L)$. One can show that
the corresponding eigenfunctions are $e^{i\l_j^+x}$ and
$e^{-i\l_j^-x}$ with
\[
\l_j^+=k\cos\theta+\frac{2\pi j}L,\quad\text{or}\quad
\l_j^-=-k\cos\theta+\frac{2\pi j}L,\quad j=0,\pm1,\pm2,\dots
\]
 We will use the system $e^{i\l_j^+x}$, which
forms an orthogonal basis in $L^2(0,L)$. One has:
\[
\int_0^L e^{i\l_j^+x}e^{-i\l_m^+x}\ dx=\int_0^Le^{\frac{2\pi i}{L}(j-m)}\
dx=0,\quad j\not=m.
\]
The normalized eigenfunctions are
\[
\p_j(x)=\frac{e^{i\l_j^+x}}{\sqrt{L}},\quad j=0,\pm1,\pm2,\dots
\]
These functions form an orthonormal basis of $L^2(0,L)$.

Let us look for $v(\bx)=v(x,y)$ of the form
\begin{equation}\label{s1_7}
v(x,y)=\sum_{j=-\infty}^{\infty} c_jv_j(y)\p_j(x), \quad y>M,\quad c_j=const.
\end{equation}
For $y>M$, equation (\ref{s1_1}) implies
\begin{equation}\label{s1_8}
v_j''+(k^2-\l_j^2)v_j=0.
\end{equation}
Let us assume that $\l_j^2\not=k^2$ for all $j$.
Then
\begin{equation}\label{s1_9}
v_j(y)=e^{i\mu_jy},
\end{equation}
where, for finitely many $j$, the set of which is denoted by $J$, one has:
\begin{equation}\label{s1_10}
\mu_j=(k^2-\l_j^2)^{1/2}>0,\quad\text{if}\quad \l_j^2<k^2, \,\, j\in J,
\end{equation}
and
\begin{equation}\label{s1_11}
\mu_j=i(\l_j^2-k^2)^{1/2},\quad\text{if}\quad \l_j^2>k^2, \,\, j\notin J.
\end{equation}

The {\bf radiation condition} at infinity requires that the scattered
field $v(x,y)$ be representable in the form (\ref{s1_7}) with $v_j(y)$
defined by (\ref{s1_9})-(\ref{s1_11}).

The {\bf Periodic Scattering Problem} consists of finding the solution to
(\ref{s1_1})-(\ref{s1_4}) satisfying the radiation condition (\ref{s1_7}),
(\ref{s1_9})-(\ref{s1_11}).

The existence and uniqueness for such a scattering problem is
established in
\cite{r461}.
 In \cite{alber} the scattering by a periodic structure was
considered earlier, and was based on a uniqueness theorem from \cite{e}.
 There are many papers
on scattering by periodic structures, of which we mention a few
\cite{alber,albertsen,bon,bonram,christiansen,kazan1,kazan2},
\\ \cite{ millar1,millar2,petit,ral1}.
The Rayleigh conjecture is discussed in several of the above
papers. It was shown (see e.g.  \cite{petit,baran})  that this
conjecture is incorrect, in general.
As we have already discussed in the previous sections,
the Modified Rayleigh Conjecture is a theorem proved in \cite{r430} for
scattering by bounded obstacles.

The main ingredient in the solution is an
analog to the
half-space Dirichlet Green's function. The function $g=g(\bx,\xi,k)$ can
be
constructed analytically ($\bx=(x_1,x_2), \bxi=(\xi_1,\xi_2)$):
\begin{equation}\label{s3_01}
g(\bx,\bxi)=\sum_j\p_j(x_1)\overline{\p_j(\xi_1)}g_j(x_2,\xi_2,k),
\end{equation}
\[
g_j:=g_j(x_2,\xi_2,k)=\begin{cases}
v_j(x_2)\psi_j(\xi_2), & \ x_2>\xi_2\\
v_j(\xi_2)\psi_j(x_2), & \ x_2<\xi_2
\end{cases}
\]
\[
\psi_j=(\mu_j)^{-1}e^{i\mu_jb}\sin[\mu_j(\xi_2+b)],\, \,
\mu_j=[k^2-\lambda_j^2]^{1/2},\quad v_j(x_2)=e^{i\mu_j x_2},
\]
where
\[
\psi_j''+(k^2-\l_j^2)\psi_j=0, \ \psi_j(-b)=0,\ W[v_j,\psi_j]=1,\
\lambda_j=k\cos (\theta)+\frac {2\pi j}{L},
\]
and $W[v,\psi]$ is the Wronskian.

The function $g$ is analytic with respect to $k$ on the complex plain with
cuts
along the rays $\lambda_j-i\tau, \, 0\leq \tau <\infty, j=0, \pm 1, \pm
2,..$, in particular, in the region $\Im k>0,$ up to the real
positive half-axis except for the set $\{\lambda_j\}_{j=0, \pm 1, \pm 2,
...}$.

{\it Choose $b>0$ such that $k^2>0$ is not an eigenvalue of the problem:}
\begin{equation}\label{s3_1}
(\Delta +k^2)\psi=0,\quad \text{in}\ D_{-b}:\,=\{(x,y)\ :\ -b\leq y\leq f(x), \quad 0\leq x\leq L\}.
\end{equation}
\begin{equation}\label{s3_2}
\begin{split}
 & \psi|_{y=-b}=0, \quad \psi_N=0\ \text{on}\ S,\\
& \psi(x+L,y)=\nu \psi(x,y),\quad \psi_x(x+L,y)=\nu \psi_x(x,y).
\end{split}
\end{equation}
One has
\begin{equation}\label{s3_3}
\begin{split}
& (\Delta +k^2)g=-\d(\bx-\bxi),\ \bx=(x_1,x_2),\ \bxi=(\xi_1,\xi_2),\\
& \bx\in \{(x,y)\ :\ -b< y<\infty, \quad 0\leq x\leq L\},
\end{split}
\end{equation}
\begin{equation}\label{s3_4}
g|_{y=-b}=0.
\end{equation}

Rayleigh conjectured \cite{ral1} ("Rayleigh hypothesis") that the series
(\ref{s1_7}) converges up to the boundary $S_L$. This conjecture is
wrong (\cite{petit}) for some $f(x)$. Since the Rayleigh
hypothesis has been widely used for numerical solution of the scattering
problem by physicists and engineers, and because these practitioners
reported high instability of the numerical solution, and there are no
error estimates, we propose a modification of the Rayleigh conjecture,
which is a Theorem. This MRC (Modified Rayleigh Conjecture) can be used
for a numerical solution of the scattering problem, and it gives an error
estimate for this solution. Our arguments are very similar to the ones
in \cite{r430}.

Rewrite the scattering problem (\ref{s1_1})-(\ref{s1_4}) as
\begin{equation}\label{s4_1}
(\Delta+k^2)v=0\ \text{in}\ D,\ v=-u_0\ \text{on}\ S_L,
\end{equation}
where $v$ satisfies (\ref{s1_4}), and $v$ has representation
(\ref{s1_7}), that is, $v$ is "outgoing", it satisfies the radiation
condition. Fix an arbitrarily small $\ep>0$, and assume that
\begin{equation}\label{s4_2}
\|u_0+\sum_{|j|\leq j(\ep)}c_j(\ep)v_j(y)\p_j(x)\|\leq\ep,
\ 0\leq x\leq L,\ y=f(x),
\end{equation}
where $\|\cdot\|=\|\cdot\|_{L^2(S_L)}$.

\begin{lemma}\label{s4_l1}
For any $\ep>0$, however small, and for any $u_0\in L^2(S_L)$, there
exists $j(\ep)$ and $c_j(\ep)$ such that (\ref{s4_2}) holds.
\end{lemma}
\begin{proof}
Let us prove
the completeness of the system
$\{\p_j(x)v_j(f(x))\}_{j=0, \pm 1, \pm2,...}$ in $L^2(S_L)$.
Assume that there is an $h\in L^2(S_L),\ h\not\equiv 0$
 such that
\begin{equation}\label{s4_3}
\int_{S_L}h\overline{\p_j(x)}v_j(f(x))\ ds=0
\end{equation}
for all $j$.
 From (\ref{s4_3}) one derives (cf. \cite{rammb1}, p.162-163)
\begin{equation}\label{s4_4}
\psi(\bx):\,=\int_{S_L}hg(\bx,\bxi) d\xi=0,\ \bx\in D_{-b}.
\end{equation}
Thus $\psi=0$ in $D_L$, and $h=\psi_N^+-\psi_N^-=0$. Lemma \ref{s4_l1}
is proved.
\end{proof}

\begin{lemma}\label{s4_l2}
If (\ref{s4_2}) holds, then
\[
\||v(\bx)-\sum_{|j|\leq j(\ep)}c_j(\ep)v_j(y)\p_j(x)\||\leq c\ep,\
\forall x,y\in D_L,
\ 0\leq x\leq L,\ y=f(x),
\]
where $c=const>0$ does not depend on $\ep, x, y$; $R>M$ is an arbitrary
fixed number, and
$\||w\||=\sup_{\bx \in D\setminus
D_{LR}}|w(\bx)|+||w||_{H^{1/2}(D_{LR})}$.
\end{lemma}
See \cite{r461} for the proof.

 From Lemma \ref{s4_l2} the basic result, Theorem \ref{s4_t1}, follows
immediately:
\begin{theorem}\label{s4_t1}
{\bf MRC-Modified Rayleigh Conjecture.} Fix $\ep>0$, however small,
and choose a positive integer $p$. Find
\begin{equation}\label{s4_9}
\min_{c_j}\|u_0+\sum_{|j|\leq p}c_j\p_j(x)v_j(y)\|:\,=m(p).
\end{equation}
Let $\{c_j(p)\}$ be the minimizer of (\ref{s4_9}). If $m(p)\leq\ep$, then
\begin{equation}\label{s4_10}
v(p)=\sum_{|j|\leq p}c_j(p)\p_j(x)v_j(y)
\end{equation}
satisfies the inequality
\begin{equation}\label{ss4_11}
\||v-v(p)\||\leq c\ep,
\end{equation}
where $c=const>0$ does not depend on $\ep$. If $m(p)>\ep$, then there
exists $j=j(\ep)>p$ such that $m(j(\ep))<\ep$. Denote
$c_j(j(\ep)):\,=c_j(\ep)$ and $v(j(\ep)):\,=v_{\ep}$. Then
\begin{equation}\label{ss4_12}
\||v-v_\ep\||\leq c\ep.
\end{equation}
\end{theorem}

\section{Numerical solution of the periodic scattering problem}
According to the MRC method (Theorem \ref{s4_t1}), if the restriction of the
incident field $-u_0(x,y)$ to $S_L$
 is approximated as in (\ref{s4_9}), then the series (\ref{s4_10})
approximates the scattered field in the entire region above the profile
$y=f(x)$. However, a numerical method that uses (\ref{s4_9}) does not
produce satisfactory results as
reported in \cite{petit} and elsewhere. Our own numerical
experiments confirm this observation. A way to overcome this difficulty
is to realize that the numerical approximation of the field
$-u_0|_{S_L}$ can be carried out by using outgoing solutions
described below.

Let $\bxi=(\xi_1,\xi_2)\in D_{-b}$, where $b>0$,
\[
D_{-b}:\,=\{(\xi_1,\xi_2)\ :\ -b\leq \xi_2\leq f(x), \quad 0\leq \xi_1\leq
L\},
\]
and $g(\bx,\xi)$ be defined as in \eqref{s3_01}. Then $g(\bx,\xi)$ is an
outgoing solution satisfying $\Delta g+k^2g=0$ in $D_L$, according to (\ref{s3_3}).

To implement the MRC method numerically one proceeds as follows:
\begin{enumerate}
\item Choose the nodes $\bx_i,\ i=1,2,...,N$
 on the profile $S_L$. These points are used to approximate
$L^2$ norms on $S_L$.
\item Choose points $\bxi^{(1)},\bxi^{(2)}, ..., \bxi^{(M)}$ in $D_{-b},\
M<N$.
\item Form the vectors ${\bf b}=(u_0(\bx_i))$, and
${\bf a}^{(m)} = (g(\bx_i,\bxi^{(m)})),
\ i=1, 2,..., N$, $m= 1,2,...,M$. Let $\bf A$ be the $N\times M$ matrix containing
vectors ${\bf a}^{(m)}$ as its columns.
\item Find the Singular Value Decomposition of ${\bf A}$. Use a predetermined $w_{min}>0$ to
eliminate its small singular values. Use the decomposition to compute
\[
r^{min}=\min\{\|{\bf b}+{\bf Ac}\|,\ {\bf c}\in \mathbb{C}^M\},
\]
where
\[
  \|{\bf a}\|^2=\frac 1{N}\sum_{i=1}^N |a_i|^2.
\]
\item {\bf Stopping criterion.} Let $\ep>0$.
\begin{enumerate}
\item If $r^{min}\leq\ep$, then stop. Use the coefficients ${\bf c}=\{c_1,c_2,...,c_M\}$
 obtained
in the above minimization step to compute the scattered field by
\[
v(x,y)=\sum_{m=1}^M c_m g(x,y,\xi^{(m)}).
\]

\item If $r^{min}>\ep$, then increase $N, M$ by the order of 2, readjust
the location of points $\xi^{(m)}\in D_{-b}$ as needed, and repeat the procedure.
 \end{enumerate}
\end{enumerate}

\begin{table}[tb]
\caption{Residuals attained in the numerical experiments on MRC for periodic structures.}
\begin{center}

\begin{tabular}{c  c  c  }

\hline
Profile & $\theta$ & $r^{min}$ \\

\hline
I & $\pi/4$ & 0.000424  \\
  & $\pi/3$ & 0.000407  \\
  & $\pi/2$ & 0.000371  \\
\hline
II & $\pi/4$ & 0.001491  \\
   & $\pi/3$ & 0.001815  \\
   & $\pi/2$ & 0.002089  \\
\hline
III & $\pi/4$ & 0.009623  \\
   & $\pi/3$ & 0.011903  \\
   & $\pi/2$ & 0.013828  \\
\hline
IV & $\pi/4$ & 0.014398  \\
   & $\pi/3$ & 0.017648  \\
   & $\pi/2$ & 0.020451  \\
\hline
\end{tabular}
\end{center}

\end{table}

We  conducted numerical experiments for four different profiles. In
each case we used $L=\pi, k=1.0$ and three values for the angle
$\theta$.  Table 5 shows the
resulting residuals $r^{min}$. Note that $\|{\bf b}\|=1$.
Thus, in all the considered cases, the MRC method
achieved $0.04\%$ to $2\%$ accuracy of the approximation.
Other parameters used in the experiments were chosen as follows:
$N=256,\ M=64,\ w_{min}=10^{-8},\ b=1.2$. The value of $b>0$, used in the
definition of $g$, was chosen experimentally, but the dependency of $r^{min}$
on $b$ was slight. The Singular Value Decomposition (SVD) is used in Step 4 since the vectors
${\bf a}^{(m)},\ m=1,2,...,M$ may be nearly linearly dependent, which
leads to an instability in the determination of the minimizer $\bf c$.
According to the SVD method this instability is eliminated by cutting
off small singular values of the matrix $\bf A$, see e.g. \cite{numrec} for details.
The cut-off value $w_{min}>0$ was chosen experimentally.
We used the truncated series (\ref{s3_01}) with $|j|\leq 120$
to compute functions $g(x,y,\xi)$. A typical run time on a 333 MHz PC was
about $40s$ for each experiment.

The following is a description of the profiles $y=f(x)$, the nodes
$\bx_i\in S_L$, and the poles $\xi^{(m)}\in D_{-b}$ used in  the computation of
$g(\bx_i,\xi^{(m)})$ in Step 3. For example, in profile I the $x$-coordinates of
the $N$ nodes $\bx_i\in S_L$ are uniformly distributed on the interval
$0\leq x\leq L$. The poles $\xi^{(m)}\in D_{-b}$ were chosen as follows:
every fourth node $\bx_i$ was moved by a fixed amount $-0.1$ parallel to
the $y$ axis, so it would be within the region $D_{-b}$. The
location of the poles was chosen experimentally to give the smallest
value of the residual $r^{min}$.

{\bf Profile I.} $f(x)=sin(2x)$ for $0\leq x\leq L,\ t_i=iL/N,\ \bx_i=(t_i,f(t_i)),\ i=1,2,...,N,\
\xi^{(m)}=(x_{4m},y_{4m}-0.1),\ m=1,2,...,M$.

{\bf Profile II.} $f(x)=sin(0.2x)$ for $0\leq x\leq L,\ t_i=iL/N,\ \bx_i=(t_i,f(t_i)),\ i=1,2,...,N,\
\xi^{(m)}=(x_{4m},y_{4m}-0.1),\ m=1,2,...,M$.

{\bf Profile III.} $f(x)=x$ for $0\leq x\leq L/2$, $f(x)=L-x$ for
$L/2\leq x\leq L,\  t_i=iL/N,\ \bx_i=(t_i,f(t_i)),\ i=1,2,...,N,\
\xi^{(m)}=(x_{4m},y_{4m}-0.1),\ m=1,2,...,M$.

{\bf Profile IV.} $f(x)=x$ for
$0\leq x\leq L,\ t_i=2iL/N,\ \bx_i=(t_i,f(t_i),\ i=1,...,N/2,
\  \bx_i=(L,f(2(i-N/2)L/N)),\ i=N/2+1,...,N,
\bxi^{(m)}=(x_{4m}-0.03,y_{4m}-0.05),\ m=1,2,...,M$.
In this profile $N/2$ nodes $\bx_i$ are uniformly distributed on its
slant part, and $N/2$ nodes are uniformly distributed on its vertical
portion $x=L$.

The experiments show that the MRC method provides a competitive
alternative to other methods for the computation of fields scattered
from periodic structures. It is fast and inexpensive. The results depend
on the number of the internal points $\xi^{(m)}$ and on their location.

\section{Inverse scattering methods based on the MRC}
Suppose that an approximate location of the obstacle $D$ is obtained a numerical inversion method,
such as the Support Function Method (SFM), see  \cite{r453,r463}. 
Then one can try to use the MRC method to improve the location of the boundary, see \cite{r430}.
Such methods are under development by the authors, and they are going to be discussed elsewhere.
Nevertheless, the MRC provides a tool for an easy construction of various examples
illustrating the severe ill-posedness of the Inverse Scattering problem,
which can be used for the algorithm's testing.

\begin{table}[tb]
\caption{Near field values of two radiating solutions with practically the same far fields.}
\begin{center}

\begin{tabular}{ r  r  r r r}

\hline
$\a'$ & $Re\ v_c$ & $Im\ v_c$ & $Re\ v$ & $Im\ v$\\

\hline
     0.00000 &  -1189.60834  &  -227.35213 &     -0.54030   &   -0.84147 \\
     0.31416 &    -73.43878  &   -15.81270 &     -0.58082   &   -0.81403 \\
     0.62832 &      1.94958  &     0.19051 &     -0.69021   &   -0.72361 \\
     0.94248 &      0.03298  &    -0.52343 &     -0.83217   &   -0.55452 \\
     1.25664 &     -1.07968  &    -0.36021 &     -0.95263   &   -0.30412 \\
     1.57080 &     -1.13445  &     0.00027 &     -1.00000   &    0.00000 \\
     1.88496 &     -0.96294  &     0.31629 &     -0.95263   &    0.30412 \\
     2.19911 &     -0.79021  &     0.55436 &     -0.83217   &    0.55452 \\
     2.51327 &     -0.66472  &     0.71819 &     -0.69021   &    0.72361 \\
     2.82743 &     -0.59154  &     0.81406 &     -0.58082   &    0.81403 \\
     3.14159 &     -0.56768  &     0.84565 &     -0.54030   &    0.84147 \\
     3.45575 &     -0.59154  &     0.81406 &     -0.58082   &    0.81403 \\
     3.76991 &     -0.66472  &     0.71819 &     -0.69021   &    0.72361 \\
     4.08407 &     -0.79021  &     0.55436 &     -0.83217   &    0.55452 \\
     4.39823 &     -0.96294  &     0.31629 &     -0.95263   &    0.30412 \\
     4.71239 &     -1.13445  &     0.00027 &     -1.00000   &    0.00000 \\
     5.02655 &     -1.07968  &    -0.36021 &     -0.95263   &   -0.30412 \\
     5.34071 &      0.03298  &    -0.52343 &     -0.83217   &   -0.55452 \\
     5.65487 &      1.94958  &     0.19051 &     -0.69021   &   -0.72361 \\
     5.96903 &    -73.43878  &   -15.81270 &     -0.58082   &   -0.81403 \\
\hline
\end{tabular}
\end{center}

\end{table}

Here we present one such example. Let the obstacle $D$ be the unit
circle $\{\bx\in\R^2\ :\ |\bx|\leq 1\}$. If the incident field is
$u_0(x)=e^{ikx\cdot\a}$, then the scattered field $v(x)=-u_0(x)$ for
$x\in S=\partial D$, and its scattering amplitude is
\be\label{s11_vfar}
A(\a',\a)=-\sqrt{\frac 2{\pi k}}\ e^{-i\frac{\pi}4}
\sum_{l=-\infty}^\infty \frac
{J_l(ka)}{H^{(1)}_l(ka)}
\  e^{il(\theta-\beta)}\,,
\end{equation}
where $\a'=\bx/\abs{\bx}=e^{i\theta}$, and $\a=e^{i\beta}$.

Let $x_1\in\R^2$. Fix an integer $L>0$, and let ${\bf c}\in \C^{2L+1}$. Form the radiating solution
\be\label{s11_vc}
v_c(x)=\sum_{l=-L}^L c_{l} H_l^{(1)}(k\abs{x-x_1})e^{il\theta_1}\,,
\end{equation}
where $(x-x_1)/\abs{x-x_1}=e^{i\theta_1}$.
 Then its far field pattern is
\be\label{s11_vcfar}
A_{v_c}(\a')=\sqrt{\frac 2{\pi k}}\ e^{-i\frac{\pi}4}
 \left(e^{-ik\a'\cdot x_1} \sum_{l=-L}^L c_{l}
\ (-i)^l e^{il\theta}\right)\,,
\end{equation}
where $\a'=x/\abs{x}=e^{i\theta}$.

Fix an $\a\in S^1$, and let
\be\label{s11_rmin}
r^{min}=\min\{\|A_{v_c}(\a')-A(\a',\a)\|\ :\ {\bf c}\in \C^{2L+1}\}\,.
\end{equation}
We conducted the minimization by the Singular Value Decomposition Method
 with the following values of the parameters:
$k=1.0,\ L=5,\ \a=(1.0, 0.0)$, and $x_1=(0.8, 0.0)$. The $L^2$ norm in
(\ref{s11_rmin}) was computed over $M=120$ directions $\a'_m$ uniformly
distributed in the unit circle $S^1$, and then normalized by $\sqrt{M}$,
so that the identity function would have the norm equal to 1. The
resulting value of the residual $r^{min}=0.00009776$ indicates that the
far field $A(\a',\a)$ was practically perfectly fit by the radiating
solution of the form (\ref{s11_vc}). However, as the Table 6 shows, the
restrictions of the exact scattered field $v$, and the fitted field
$v_c$ to the boundary $S$ of the obstacle $D$ are vastly different. The
columns in Table 6 correspond to the real and the imaginary parts of the
scattered fields, and the rows correspond to different values of the
angle $\a'$. Thus, one has to conclude that, as expected, a coincidence
of the radiating solutions at the far field does not imply that the near
 fields are also coincidental, see \cite{r463}.

\section{Stability Index Method}\label{section_sim}

Various algorithms for direct and inverse scattering problems require
global minimization of functions of many variables, see \cite{r463}.
Since most objective functions contain many local minima, this is a
highly nontrivial task. In several papers, starting with [14], the authors 
developed and
tested the Stability Index Method (SIM) for global minimization. In our
presentation here we follow \cite{g8}, which also  
contains a convergence analysis 
and additional numerical results.

The Stability Index Method  combines stochastic and
deterministic algorithms to find global minima of multidimensional
functions. The functions may be nonsmooth and may have multiple
local minima. The method examines the change of the diameters of
the minimizing sets for its stopping criterion. At first,
the algorithm uses the uniform random distribution in the admissible
set. Then normal random distributions of decreasing variation are
used to focus on probable global minimizers. To test the method,
we have applied it to standard test functions of several
variables. The computational results show that the SIM is
efficient, reliable and robust.

Given a function $f : A\ra \R$, our goal is to minimize it over an admissible
set $A$ assumed to be a bounded set in a metric space $X$.
Typically, the
structure of the objective function $f$ is quite complicated. In particular, it can
 have
many local minima and a non unique global minimum. To better understand the
structure of the minima, let us introduce the minimizing sets $S_\ep$ of $f$.
Let $m=\inf\{f(x) : x\in A\}$. Given an $\ep>0$ define
\begin{equation}\label{i1_e1}
S_{\ep}=\{x\in A\ :\ f(x)<m+\ep \},
\end{equation}
or
\begin{equation}\label{i1_e2}
S_{\ep}=\{x\in A\ :\ f(x)<f(x_p)+\ep \},
\end{equation}
if the problem admits a global minimizer $x_p\in A$.

{\bf Definition}. {\it Given an $\ep>0$, let $D_\ep$ be the diameter of the minimizing set
$S_\ep$. We call  $D_\ep$ the {\bf Stability Index} of the
minimization
problem (\ref{i1_e1}).}

We are interested in the
behavior of $D_\ep$  as $\ep\ra 0$. So, one can  say that the
problem \eqref{i1_e1} possesses a set of Stability Indices $\{D_\ep :
\ep>0\}$, and the above definition should be understood in this sense.

One would expect to obtain a stable identification for minimization
problems with small (relative to the admissible set) stability indices. Minimization
problems with large stability indices either have distinct global minimizers, or the function $f$
is nearly flat in a neighborhood of the global minimizer $x_p$.
In this situation, and with no additional
information known, one has an uncertainty
of the minimizer's choice. The stability index provides a quantitative
measure of this uncertainty or instability of the minimization.

In a practical minimization problem one constructs a sequence of minimizers
$\{x_1,x_2,...\}\subset A$, and makes a decision when to
terminate the iterations according to a stopping criterion. We assert
that the knowledge of the Stability Index provides a valuable tool for
the formulation of such a stopping criterion.

Originally, we have applied the Stability Index minimization method to
inverse scattering problems arising in quantum mechanical
scattering, \cite{r435}.
Such potential scattering
problems are important in quantum mechanics, where they
appear in the context of scattering of particles bombarding
an atom nucleus. One is interested in reconstructing the
scattering potential from the results of a scattering
experiment. Assuming a particular structure of the potential, the
scattering results can be computed and compared to the given scattering
data. Thus the inverse scattering problem is reduced to the minimization
of the discrepancy (best fit to data), see \cite{r435,r479} for details.

The goal of the SIM algorithm is to find a minimizing set $S_\ep$ that
fits within a small portion of the computational domain $A\subset \R^N$. Practically, we assume
that
 $A = [-M,M]^N\subset \R^N$, for an
$M>0$. If it is desirable to introduce different scales for the variables, then the algorithm
 should be modified accordingly.

Let $0<\delta<1$. The minimization is stable if, given a global minimizer $x_p$,
 we are able to find a minimizing set
$S_\ep\subset C[x_p,\delta]$, where $C[x_p,\delta]$ is the cube centered
at $x_p\in A$ with the side equal to $2\delta M$.

The next step is to define a sequence of normal distributions $T_n$
with the variances $\mu_n\ra 0$, as $n\ra\infty$. Thus we fix an
$0<\alpha<1$, and let $\mu_n=\alpha^n,\; n=1,2,...$

Initially, for $n=0$, let the $T_0$ be the uniform random distribution    in $A$.
A special algorithm {\bf SMS}, described below,  determines a finite Stable
Minimizing Set $S_0\subset A$. Let $x_0\in S_0$ be the minimizer in $S_0$, that is
\begin{equation}
f(x_0)=\min\{f(x) : x\in S_0\}\,.
\end{equation}
If $S_0\subset C[x_0,\delta]$, then the minimization is stable
and the global minimizer $x_p=x_0$.

If, on the other hand, the above inclusion is not achieved, then one
continues with another application of the {\bf SMS}, this time using the
normal distribution $T_1$ with the mean at $x_0$, and the variance
$\mu_1$, etc. The iterations continue until either $S_n\subset
C[x_n,\delta]$ or $3\mu_n<2\delta M$. The last condition is needed to prevent
all the trial points to be chosen too close to $x_n$, thus preventing a
reasonable estimate for the diameter of $S_n$.

\vskip2mm
\leftline{\bf Stability Index Method (SIM)}

Fix $0<\a,\delta<1$. Suppose that
 $A= [-M,M]^N $.

\begin{enumerate}
\item  Initialization. Let $n=0$. Use the {\bf SMS} algorithm with
 the uniform random distribution $T_0$
 in $A$ to determine the minimizing set $S_0\subset A$ and the minimizer
$x_0\in S_0$.
Go to the Stopping Criterion (step 3) to check if additional iterations are needed.

\item ($n-$th iteration). Let $\mu_n=\alpha^n$. Use the {\bf SMS} algorithm with
the normal random distribution $T_n$ with the mean at $x_{n-1}$ and the variance
$\mu_n$
 to determine the minimizing set $S_n\subset A$ and the minimizer
$x_n\in S_n$.

\item Stopping criterion.
Let $C[x_n,\delta]$ be the cube centered
at $x_n\in A$ with the side equal to $2\delta M$.

If $S_n\subset
C[x_n,\delta]$, then stop. The minimization is stable.
The estimated global minimizer $x_p$ is $x_n$.

If $S_n\not\subset
C[x_n,\delta]$ and $3\mu_n<2\delta M$, then stop. The minimization is unstable.
The diameter (Stability Index) $D_n$ of $S_n$ is a
measure of the instability of the minimization.

Otherwise,  increase $n$ by $1$, and return to Step 2
to do another iteration.

\end{enumerate}

Note that the obtained point $x_p$ is an estimated global minimizer.

The main part of the Stability Index Method  is the {\bf SMS} algorithm
which determines stable minimizing sets $S_n$, corresponding to the
random distributions $T_n$. These distributions are either uniform in
$A$ or normal with a given variation $\mu_n$.

The {\bf SMS} algorithm is  an iterative algorithm.
It can be called
 an Iterative Reduced Random Search method. Choose an integer $K>0$ 
assuming that $K$ random points in $S_\ep$ are sufficient
to  estimate  its diameter $D_\ep$. If $n\geq 1$, then the calling
algorithm SIM provides the minimizing set $S_{n-1}$, its minimizer
$x_{n-1}$, and the variance $\mu_n$.

Let a batch $H^1\subset A$ of $L>K$ trial points be generated in the admissible set
$A$ according to the random distribution $T_n$. If $n=0$, then $T_0$ is
just the uniform random distribution in $A$. If $n\geq 1$, then $T_n$ is the
normal distribution with the variance $\mu_n$, and the mean at $x_{n-
1}$.
Let $Q^1_U$ be the subset of $K$ points from $H^1$ where the objective
function $f$ attains its $K$ smallest values. That is
\begin{equation}
\max\{f(u_i) : u_i\in Q^1_U\}\leq \min\{f(u_i) : u_i\not\in Q^1_U\}.
\end{equation}

Use each point $u_i\in Q^1_U$ as the initial guess for a Local
Minimization Method ({\bf LMM}) of your choice, e.g. the conjugate gradient method, etc.
 The specific {\bf LMM}
used by us is described below. While the use of a local
minimization is not, strictly speaking, necessary for the SIM, it provides a
significant improvement in the performance of the algorithm, and is
highly recommended. Thus for each starting point $u_i\in Q^1_U$ the {\bf LMM}
 produces a minimizer
$v_i\in A$.
Let $Q^1_V$ be the set of all such minimizers.
Let $Q^1$ be the subset of $Q^1_U\cup Q^1_V$ containing $K$ points with
the smallest values of $f$, and $q^1$ be the minimizer in $Q^1$. Define
the radius of $Q^1$ by

\be
R^{(1)}=\max\{\|z_i-q^1\| : z_i\in Q^1,\; i=1,2,...,K\}.
\end{equation}

The
idea of the Stability Index Method is to iteratively construct subsets
$Q^j$ until their diameters are stabilized. Practically, one can achieve
the same goal by estimating and examining the radius $R^{(j)}$ of the set $Q^j$.
This also requires less computational effort.

To construct the next set $Q^2$ generate another batch $H^2\subset A$ of $L$ trial points
according to the uniform random distribution, if $n=0$, or, for $n\geq 1$,
according to the normal distribution $T_n$ with the variance $\mu_n$,
and the mean at $q^1$. Let $Q^2_U$ be the subset of $K$ points from
$H^2\cup Q^1$ having the smallest $K$ values of $f$. Apply the {\bf LMM} to produce the set of minimizers
$Q^2_V$. Of course, if some point $u_i\in Q^2_U$ has already been used
as an initial guess for the {\bf LMM} in the previous iteration, it is
excluded from the {\bf LMM} application.  Let $Q^2$ be the subset of $Q^2_U\cup Q^2_V$ containing $K$ points with
the smallest values of $f$. Let $q^2$ be the minimizer in $Q^2$, and
$R^{(2)}=\max\{\|z_i-q^2\| : z_i\in Q^2,\; i=1,2,...,K\}$ be its
radius, etc.

This way one produces a sequence of the minimizing sets $Q^j,\;
j=1,2,...$. Let $0<\gamma<1$, and $P$ be a positive integer.
 The iterations are terminated if the maximum number
of iterations $N_{max}$ is exceeded or the following Stopping Criterion
is satisfied:

\be\label{i2_stop}
\left|R^{(j)}-\frac 1P\sum_{i=j-P+1}^j R^{(i)}\right|<2\gamma M\,.
\end{equation}

In either case, when the last iteration $j$ is determined from
\eqref{i2_stop} or $j=N_{max}$, we let $S_n=Q^j$ and $x_n=q^j$.

\vskip2mm

\leftline{\bf Stable Minimizing Set (SMS) algorithm}

Fix $0<\gamma<1$, and integer $K,L>K,P,N_{max}$. Constant $M$, normal random distribution $T_n$,
  its variance
$\mu_n$ (for $n\geq 1$), the minimizing set $S_{n-1}$, and the minimizer $x_{n-1}$
 are supplied by the calling algorithm SIM.

\begin{enumerate}
\item Initialization. Let $j=1$.

\begin{itemize}

\item For $n=0$.  Generate a batch $H^1$ of $L$
trial points in $A\subset\R^N$ using the
uniform random distribution. Let $Q^1_U$ be the subset of $K$ points from $H^1$ where the objective
function $f$ attains its $K$ smallest values. Go to step 4.

\item For $n\geq 1$. Generate a batch $H^1$ of
$L$ trial points in $A\subset\R^N$ using the normal distribution $T_n$
with the variance $\mu_n$ and the mean at $x_{n-1}$.
 Let $Q^1_U$ be the subset of $K$ points from $H^1\cup S_{n-1}$ where the objective
function $f$ attains its $K$ smallest values. Go to step 4.
\end{itemize}

\item Iterative step ($j\geq 2$).

\begin{itemize}

\item For $n=0$.  Generate a batch $H^j$ of $L$
trial points in $A\subset\R^N$ using the
uniform random distribution.

\item For $n\geq 1$. Generate a batch $H^j$ of
$L$ trial points in $A\subset\R^N$ using the normal distribution $T_n$
with the variance $\mu_n$ and the mean at $q^{j-1}$.
\end{itemize}

\item Let $Q^j_U$ be the subset of $K$ points from $H^j\cup Q^{j-1}$ where the objective
function $f$ attains its $K$ smallest values.

\item Local minimization. Use each unflagged point $u_i\in Q^j_U$ as the initial guess for a Local
Minimization Method ({\bf LMM}). Let $v_i\in A$ be the resulting
minimizer.
Let $Q^j_V$ be the set of all such minimizers resulting from the application of {\bf LMM} to
$Q^j_U$. Flag all points in $Q^j_U$ and $Q^j_V$.

\item
Let $Q^j$ be the subset of $Q^j_U\cup Q^j_V$ containing $K$ points with
the smallest values of $f$ and $q^j$ be the minimizer in $Q^j$. Define
the radius of $Q^j$ by
\[
R^{(j)}=\max\{\|z_i-q^j\| : z_i\in Q^1,\; i=1,2,...,K\}.
\]

\item Stopping criterion.
\begin{itemize}

\item If $j<P$, increase $j$ by $1$ and return to step 2 for another
iteration.

\item If $j\geq P$, compute the average radius during the last $P$ iterations:
\[
R_a=\frac 1P\sum_{i=j-P+1}^j R^{(i)}.
\]

\item Termination. If $|R^{(j)}-R_a|\leq 2\gamma M$, or $j\geq N_{max}$,  let $S_n=Q^j$, $x_n=q^j$
and exit the procedure.

\item Otherwise, increase $j$ by $1$ and return to step 2 for another
iteration.
\end{itemize}
\end{enumerate}

The {\bf SMS} implementation involves a combination of stochastic
(global)
and deterministic (local) minimization methods. 
Generally, local searches offer more precision
and speed than their global counterparts, so that adding a local step to
a global minimization algorithm should yield improvement in both areas.
Likewise, by itself, a local minimization method will very often produce
points of considerable distance from the actual global minimizer, that is
it would be trapped in one of many local minima of the objective function $f$.
 Adding
a global step  helps the algorithm escape from local minima, and
explore the entire admissible set $A$. The use of various normal
distributions of decreasing variance is similar to ideas of the
simulated annealing method \cite{kir}.

The particular Local Minimization Method (LMM)  used in the numerical
experiments was a modification
of Powell's minimization method in $\R^N$, \cite{bre}. It was chosen
with applications in mind, for which the objective function $f$ does not have a convenient
expression for its gradient. Either a Golden Search or Brent's method
can be used for one dimensional minimizations, \cite{numrec}.

\subsection*{Modified Powell's Method}
\begin{enumerate}
\item  Choose the set of directions $u_i\,,\;i=1,2,\dots,N,$ to be the 
standard basis
 in $\R^N$
\[
u_i=(0,0,\dots,1,\dots,0)\,,
\]
where $1$ is in the i-th place.
\item  Save the starting point $p_0$.
\item  For $i=1,\dots,N$ move from
$p_{i-1}$ along the direction $u_i$ and find the point
of minimum $p_i$.
\item  Set $v=p_{N}-p_0$.
\item  Move from $p_0$  along the direction $v$ and find the minimum. Call it
 $p_0$ again. It replaces $p_0$ from step 2.
\item Repeat the above steps until a stopping criterion is satisfied.
The resulting point is $p_{min}$.
\end{enumerate}
Note that $f(p_{min})\leq f(p_0)$ for any objective function $f$ used in the
Local Minimization Method.

\section{Numerical results for SIM}

    The Stability Index Method described in the previous sections
was tested on several functions designed to test and compare various  minimization algorithms, see \cite{g8}
for additional test functions results.
The experiments were conducted on a 2.8 GHz PC with 256 MB RAM.

In all the numerical experiments we used the same parameter values: $\alpha=0.8, \delta=0.001,
\gamma=0.001, K=30, L=5000, P=6$, and $N_{max}=30$.
For each test function the admissible set $A$ is a
cube $[-M,M]^N$ provided in the function's description
together with its global minimizer.

\vskip5mm
\leftline{\it Test Function 1}
\begin{multline*}
f(x,y)=\left(\sum_{i=1}^5 i\cos[(i+1)x+i]\right)
\left(\sum_{i=1}^5 i\cos[(i+1)y+i]\right)\\
+0.5((x+1.4213)^2+(y+0.80032)^2)
\end{multline*}
The minimum is sought on $[-5,5]\times[-5,5]$. This function has
a global minimum at $(-1.42513, -0.80032)$ with a function value
of $-186.73091$, \cite{yiu}.

\vskip5mm
\leftline{\it Test Function 2}
\begin{multline*}
f(x,y)=e^{\sin(50x)}+\sin(60e^y)+\sin(70\sin x)+\sin(\sin(80y))\\
-\sin(10(x+y))+(x^2+y^2)/4.
\end{multline*}
   The minimum is sought on $[-1,1]^2$. According to \cite{bornemann} the minimum occurs at approximately
$(-0.0244031, 0.2106124)$ with a function value of $-3.30686865$.

\vskip5mm
\leftline{\it Test Function 3}
\begin{equation*}
f(x)=\frac{\pi}{N}\left(10\sin^2(\pi y_1) +
\sum_{i=1}^{N-1}((y_i-1)^2(1+10\sin^2(\pi y_i+1))+(y_N-1)^2\right)
\end{equation*}
where $x=(x_1,x_2,...,x_N)\in\R^N,\; y_i=1+0.25(x_i-1),\; i=1,2,...,N$.
The minimum is sought on $[-10,10]^N$. This function has a
global minimum at $x=(1,1,...,1)$ with a function value
of $0$, \cite{yiu}.

The results of the minimization using the {\bf SIM} for
these test functions are shown in Table 7.
   The algorithm was run 20 times for each function.
It found the correct global minimum most of the time.
The "success rate" column in Table 7 shows the percentage of trials
in which the global minimum was found exactly. The "Function evaluation"
column shows the average number of times the objective function was
evaluated.
 Finally, Table 7 shows the average run
time, in seconds, for a single trial run.

\begin{table}[tb]
\caption{Results of the computational experiments for SIM.}
\begin{center}

\begin{tabular}{c  c  r  r  c   c}

\hline
Function&  Dimension &  Actual   &   Found &       Success &  Average run \\
    	&      $N$   &  minimum  &  minimum&     rate (\%) &  time (seconds)\\
\hline
1&      2&   -186.731&   -186.731&  100&    2\\
2&      2&   -3.30687&   -3.30687&  100&    2\\
3&      5&    0.00000&  0.00000&    100&    7\\
3&      10&   0.00000&  0.00000&    100&    16\\
3&      20&   0.00000&  0.00000&    100&    50\\
\hline
\end{tabular}
\end{center}

\end{table}

\section{Conclusions}
Let $D$ be a 2D or 3D obstacle, $S$ be its boundary, and $u_0$ be the incident field. 
Rayleigh conjectured that the acoustic field
$u$ in the exterior of the obstacle is given by 
\be
u(x, \alpha) = e^{ik\alpha \cdot x} + 
\sum^\infty_{\l =0} A_\l (\alpha)\psi_\l, \quad
\psi_\l:= Y_\l (\alpha^\prime) h_\l (kr), \quad  
\alpha^\prime =
  \frac{x}{r}.  
\end{equation}
 While this conjecture (RC) is
false for many obstacles, it has been modified to obtain a
representation for the solution of
(\ref{s1_e1})-(\ref{s1_e2}) and to obtain its error.

It is proved  that if
 $v_\ep$ is an outgoing solution of the Helmholtz equation in the exterior 
domain  $D'$ and  
 $u_0 +v_\ep$ approximates zero in $L^2(S)-$norm on the boundary $S$,  
then $v_\ep$   approximates the exact scattered
field $v$ in  $D'$, see Theorems \ref{s2_thm1} and \ref{s2_thm2}.
The Modified Rayleigh Conjecture approach to obstacle scattering
problems is based on the following observation: 
the functions $\psi_\l(x,z),\;z\in D$
and their linear combinations are outgoing solutions to
the Helmholtz equation in the exterior domain. Therefore, one
just needs to find a  combination of such functions that gives the best fit to 
$-u_0$ on the boundary $S$. Then this combination 
approximates the scattered field everywhere in the exterior $D'$ of
the obstacle $D$ and the error of this approximation is given in 
Theorems \ref{s2_thm1} and \ref{s2_thm2}.

In this paper we describe several  implementations of the MRC method
which give an efficient approach to solving obstacle scattering problems
for 2D and 3D problems with complicated geometries.
Our implementations of the MRC method worked  more 
efficiently than the BIEM method.

Various methods for solution of direct and inverse scattering problems
require a global minimization of the objective function. We developed 
the Stability Index Method which is a robust and efficient algorithm for
global minimization. Its efficiency comes from a combined use of
global and local minimization. The global (stochastic) part
employs uniform and normal random distributions. It can be
combined with local (deterministic) methods appropriate for the
objective function. The diameters of the minimizing sets
(Stability Index) are used for a self-contained stopping
criterion. The computational experiments show that the method was
successful for various standard test functions over
multidimensional domains. 
No adjustment of parameters was needed
in different tests. The method is well suited for low dimensional
minimization problems. Its performance  deteriorates for higher dimensional problems.
 The Stability Index Method is a valuable
addition to already existing global minimization methods.

\end{document}